\documentclass[10pt]{amsart}
\usepackage[USenglish]{babel}
\usepackage{amsmath, amsfonts, amssymb, amsthm, amscd}
\usepackage{bbm, caption, subcaption, color, csquotes, dsfont, enumitem, epsf, graphicx, indentfirst, latexsym, mathtools, rotating, scalerel, tikz, tikz-cd, verbatim}
\usepackage[final, breaklinks=true, colorlinks, citecolor=green, linkcolor=blue, pdfusetitle]{hyperref}

\usetikzlibrary{decorations.pathmorphing}

\usepackage[T1]{fontenc}

\usepackage[abbrev]{amsrefs}

\usepackage{lmodern}

\newcommand{\bm}[1]{\boldsymbol{#1}}

\newtheorem{theorem}{Theorem}[section]

\newtheorem{proposition}[theorem]{Proposition}
\newtheorem{corollary}[theorem]{Corollary}

\newtheorem*{mainresult}{Main result}

\newcommand{\thisaxiomname}{}
\newtheorem*{genericaxiom*}{\thisaxiomname}
\newenvironment{axiom*}[1]{\renewcommand{\thisaxiomname}{#1}%
	\begin{genericaxiom*}}{\end{genericaxiom*}}

\theoremstyle{definition}
\newtheorem{definition}[theorem]{Definition}

\newtheorem{problem}{Problem}
\newtheorem{solution}{Solution}

\theoremstyle{remark}

\numberwithin{equation}{section}

\newcommand{\R}{\mathbb{R}}
\newcommand{\Z}{{\mathbb Z}}

\newcommand{\C}{{\mathbb C}}

\newcommand{\Q}{{\mathbb Q}}

\newcommand{\cA}{{\mathcal A}}
\newcommand{\cB}{{\mathcal B}}

\newcommand{\cE}{{\mathcal E}}

\newcommand{\cM}{{\mathcal M}}

\newcommand{\cS}{{\mathcal S}}

\newcommand{\cU}{{\mathcal U}}

\newcommand{\cW}{{\mathcal W}}
\newcommand{\cX}{{\mathcal X}}

\newcommand{\cZ}{{\mathcal Z}}

\newcommand{\CM}{\overline{\mathcal{M}}}

\newcommand{\ww}{\omega}

\newcommand{\delbar}{\overline{\partial}}
\newcommand{\delbarj}{\overline{\partial}_J}

\newcommand{\delbarjinv}{\smash{\overline{\partial}_J}\vphantom{\partial}^{-1}}

\DeclareMathOperator{\ind}{ind}

\DeclareMathOperator*{\GW}{GW}
\DeclareMathOperator{\PD}{PD}

\newcommand{\ssc}{\text{sc}}
\newcommand{\dmspace}{\overline{\mathcal{M}}}
\newcommand{\dmlog}{\smash{\overline{\mathcal{M}}}\vphantom{\mathcal{M}}^{\text{log}}}

\newcommand{\id}{\operatorname{id}}

\newcommand{\sign}{\text{sgn}}

\title[Pseudocycle are strict subset of polyfold {GW}-invariants]{Pseudocycle {G}romov--{W}itten invariants are a strict subset of polyfold {G}romov--{W}itten invariants}
\date{\today}
\thanks{Research partially supported by Project C5 of SFB/TRR 191 ``Symplectic Structures in Geometry, Algebra and Dynamics,'' funded by the DFG}

\author{Wolfgang Schmaltz}
\address{Faculty of Mathematics, Ruhr-Universit{\"a}t Bochum, 44801 Bochum, Germany}
\email{\href{mailto:wolfgang.schmaltz@rub.de}{wolfgang.schmaltz@rub.de}}
\urladdr{\url{https://sites.google.com/view/wolfgang-schmaltz/home}}

\begin{document}


\begin{abstract}
	Given a semipositive symplectic manifold, we prove that the pseudocycle genus-zero Gromov--Witten invariants are equal to the polyfold genus-zero Gromov--Witten invariants.
\end{abstract}

\maketitle

\tableofcontents


\section{Introduction}

In the decades since Gromov's seminal paper \cite{G}, Gromov--Witten invariants have remained a topic of lasting fascination in the symplectic geometry community.

The basic geometric idea of a GW-invariant is simple enough to grasp: given a symplectic manifold, it is the finite count of the $J$-holomorphic curves which, at a given marked point, pass through a specified submanifold of the symplectic manifold.
Due to the insights of Gromov, the space of $J$-curves can be given the structure of a compact topological space by adding additional nodal curves.
If we imagine, for a moment, that the relevant space of $J$-curves also has the structure of a \textit{compact manifold}, this finite count can be interpreted as an intersection number between the image of this compact manifold and the submanifolds of our symplectic manifold.
Thus the problem becomes: \textit{How can we give the space of $J$-curves the structure of something like a compact manifold, in particular, a space with sufficient structure for defining intersection numbers?}

If the Cauchy--Riemann section were completely transversal to the zero section, i.e., transversal when considered on all possible nodal strata of the space of $J$-curves,
then the top stratum of the space of $J$-curves would have the structure of a finite-dimensional manifold and all nodal strata would have the structure of manifolds with codimension at least $2$ relative to the top stratum.
Thus the space of $J$-curves would have the structure of a ``pseudocycle'', i.e., a space whose boundary (the nodal strata) would be \textit{invisible} from the point of view of homology.
Moreover, pseudocycles have sufficient structure for defining intersection numbers, as for dimension reasons the boundary will not contribute to the intersection number.

However, transversality of the Cauchy--Riemann section is often impossible to obtain through classical techniques, i.e., perturbation of an almost complex structure.
As a consequence, the nodal strata may have dimension larger than expected---indeed, larger than the dimension of the top stratum.
This situation is a fundamental obstacle to the rigorous definition of a GW-invariant; after all, for dimension reasons, any (nonzero) contribution to the intersection number from such a nodal stratum with large dimension would no longer be finite.

The pseudocycle approach to defining GW-invariants deals with this difficulty by imposing strict conditions on the symplectic manifold which disallows such phenomena in the nodal strata.
The ``semipositive'' condition was first introduced by McDuff in 1991 in \cite{mcduff1991symplectic}, and is specifically designed to guarantee that in the genus-zero case, the strata of nodal $J$-curves will have codimension at least $2$ relative to the dimension of the top stratum of non-noded simple $J$-curves.
Thus for semipositive symplectic manifolds, the space of genus-zero $J$-curves has the structure of a pseudocycle.
Around the end of 1993, GW-invariants for genus-zero were first rigorously defined for semipositive symplectic manifolds by Ruan--Tian in \cite{rt1995quatumcohomology}.
A complete description of the theory of pseudocycle GW-invariants is contained in \cite{MSbook}, and a detailed treatment of the theory of pseudocycles independent of symplectic geometry may be found in \cite{Z}.

The polyfold approach to defining GW-invariants is to achieve transversality of the Cauchy--Riemann section directly by means of an ``abstract perturbation''. This approach does not require any hypothesis on the symplectic manifold nor on the genus.
The resulting zero set of a perturbed Cauchy--Riemann section has a global smooth structure of a compact oriented ``weighted branched orbifold'', which is sufficient structure to define invariants.
Polyfold theory was used by Hofer--Wysocki--Zehnder in \cite{HWZGW} to give a well-defined GW-invariant for arbitrary genus and for general symplectic manifolds via the integration of differential forms. 
This approach was extended by Schmaltz in \cite{schmaltz2019steenrod} to equivalently define the polyfold GW-invariants via intersection numbers.

In this paper, we show that these two approaches yield GW-invariants which are equal.

\begin{mainresult}
	For a given semipositive symplectic manifold, the pseudocycle genus-zero Gromov--Witten invariants are equal to the polyfold genus-zero Gromov--Witten invariants.
\end{mainresult}

Since the polyfold GW-invariants are also defined for general symplectic manifolds and for arbitrary genus, we have
\[
	\left\{ \begin{array}{c} \textit{pseudocycle} \\ \textit{Gromov--Witten} \\ \textit{invariants} \end{array} \right\}
	\subsetneq
	\left\{ \begin{array}{c} \textit{polyfold} \\ \textit{Gromov--Witten} \\ \textit{invariants} \end{array} \right\}.
\]

In \S\ref{sec:GW-invariants} we discuss the pseudocycle and polyfold approaches to GW-invariants in greater detail.
In \S\ref{sec:the_abstract_comparison_thm} we prove an abstract comparison theorem which makes it possible to compare intersection numbers for pseudocycles with intersection numbers for compact oriented weighted branched orbifolds, and moreover allows us to show that the pseudocycle GW-invariants are equal to the polyfold GW-invariants in some situations.
In \S\ref{sec:pseudocycle-full-gw-invariant} we discuss the details of the complete definition of the pseudocycle GW-invariants which requires us to consider the space of $J$-curves of a graph as well as $J$-curves modeled on fixed trees.
In \S\ref{sec:polyfold_of_graphs} we show how the setup of $J$-curves of graphs modeled on fixed trees can be interpreted within polyfold theory.
In \S\ref{sec:equality_polyfold_GW_invariants} we prove that the variously defined polyfold GW-invariants are equal, thus completing the proof of our main result.
Finally in \S\ref{sec:the_mixed_GW-invariants_of_RT} we show that the mixed Gromov--Witten invariants of Ruan--Tian are equal to the correspondingly defined polyfold Gromov--Witten invariants.


\section{Gromov--Witten invariants}
\label{sec:GW-invariants}

\subsection{What is a Gromov--Witten invariant?}

Let $(M,\ww)$ be a closed symplectic manifold of dimension $\dim M = 2n$, and let $J$ be a $\ww$-compatible almost complex structure.
For a spherical homology class $A\in H_2(M;\Z)$, for genus $g=0$, and for an integer $k\geq 3$, we consider the set of $J$-holomorphic curves:
\[
	\cM_{A,0,k} (J) :=
	\left\{
	\begin{array}{c}
		u: S^2 \to M \\
		\{z_1,\ldots,z_k\} \in S^2
	\end{array}
	\middle|
	\begin{array}{c}
		\tfrac{1}{2} (du + J \circ du \circ i) = 0 \\
		u_* [S^2] = A, z_i \neq z_j \text{ if } i \neq j
	\end{array}
	\right\}
	\biggm/
	\text{PSL}(2,\C).
\]
Gromov proved this set has a natural compactification in \cite{G}, which was later refined into the stable map compactification of Kontsevich in \cite{Kstable}, and thus we may also consider its \text{compactification}, the \textbf{Gromov--Witten moduli space}:
\[
	\CM_{A,0,k} (J) := \cM_{A,0,k} (J) \sqcup \{\text{stable nodal $J$-holomorphic curves}\}.
\]
We seek to use this space to construct invariants of the symplectic manifold $M$.
To this end, we define the \textbf{evaluation map} which evaluates a stable curve on each marked point:
\[
	ev: \CM_{A,0, k} (J) \to M\times \cdots \times M.
\]
On the top stratum of non-noded stable $J$-holomorphic curves it is given by
\[
	ev\left([(u,z_1,\ldots,z_k)] \right): = (u(z_1),\ldots, u(z_k)).
\]

We also consider the \textbf{Grothendieck--Knudsen space}, the natural compactification of the space of configurations of $k$-marked points on a sphere modulo biholomorphic equivalence:
\[
	\dmspace_{0,k} := \overline{\{ \{z_1,\ldots,z_k\}\in S^2 \mid z_i\neq z_j \text{ if } i \neq j\} / \text{PSL}(2;\C) }.
\]
It is the special case for genus $g=0$ of the more general Deligne--Mumford spaces $\dmspace_{g,k}$.
In contrast to the Deligne--Mumford spaces (which are finite-dimensional oribfolds), the Grothendieck--Knudsen space is a finite-dimensional manifold of dimension $\dim \dmspace_{0,k} = 2k-6$.
We may define a \textbf{projection map} from the GW-moduli space to the GK-space which forgets the curve which maps to $M$ and which stabilizes the resulting unstable domain components:
\[
	\pi: \CM_{A,0,k} (J) \to \dmspace_{0,k}.
\]
On the top stratum of non-noded stable $J$-holomorphic curves it forgets the map $u$ and is given by
\[
	\pi\left([(u,z_1,\ldots,z_k)]\right) := [(z_1,\ldots,z_k)].
\]

The traditional interpretation of a \textbf{Gromov--Witten invariant} is the (supposedly) finite count of $J$-holomorphic curves which at the $i$th-marked point pass through a submanifold $\cX_i \subset M$ and whose marked point configuration is restricted by the projection map to a submanifold $\cB \subset \dmspace_{0,k}$.

Such an intersection number should depend only on the homology classes of the submanifolds, and should be independent of the almost complex structure. This count can be packaged algebraically as a homomorphism:
\[
	\GW_{A,0,k} : H_*(M;\Q)^{\otimes k} \times H_*(\dmspace_{0,k};\Q) \to \Q.
\]

A foundational problem in symplectic geometry is to actually show that such a GW-invariant is \textit{well-defined}.
Ideally, we would like to define a GW-invariant rigorously via an intersection number:
\[
	\GW_{A,0,k} ([\cX_1],\ldots,[\cX_k];[\cB]) = ``(ev\times \pi) (\CM_{A,0,k} (J)) \cdot (\cX_1\times \cdots \times \cX_k \times \cB) ",
\]
or as an integral:
\[
	\GW_{A,0,k} ([\cX_1],\ldots,[\cX_k];[\cB]) = ``\int_{\CM_{A,0,k}(J)} ev^* (\PD [\cX_1]\wedge \cdots \wedge \PD [\cX_k]) \wedge \pi^* \PD [\cB]",
\]
or as a pairing with a (virtual) fundamental class:
\[
	\GW_{A,0,k} ([\cX_1],\ldots,[\cX_k];[\cB]) = ``\left\langle  (ev\times\pi)_* [\CM_{A,0,k}(J)], \PD [\cX_1\times\cdots\times\cX_k\times\cB] \right\rangle".
\]
Such definitions require additional structure on the GW-moduli space; an intersection number requires tangent spaces and notions of transversal intersection, an integral requires smooth partitions of unity and notions of differential forms, and a (virtual) fundamental class requires a distinguished homology class on the topological space.

However, \textit{a priori}, the GW-moduli space only has the structure of a compact topological space, and this alone is insufficient to define any of the above. More structure is needed.


\subsection{Gromov--Witten invariants via pseudocycles}

A symplectic manifold $(M^{2n},\ww)$ is called \textbf{semipositive} if, for every $A\in \pi_2(M)$,
\begin{equation} \label{eq:semipositive}
	\ww(A) >0,\ c_1(A) \geq 3-n \quad \implies \quad c_1(A) \geq 0.
\end{equation}
As we will see, the semipositive condition allows for precise control on the dimension of the nodal $J$-curve strata of the GW-moduli spaces. This precise control allows us to give the GW-moduli spaces the structure of a \textit{pseudocycle}.

\subsubsection{Geometric approach to regularization}

To begin, we recall the so-called ``geometric approach'' to regularization of the GW-moduli space. This is the approach described in full detail in \cite{MSbook}.

Consider the Banach manifold
\[
	\cB^{m,q}_{A,0,k} :=
	\left\{
	\begin{array}{c}
		u: S^2 \to M \\
		\{z_1,\ldots,z_k\} \in S^2
	\end{array}
	\middle|
	\begin{array}{c}
		u\in W^{m,q},\ u_* [S^2] = A \\
		z_i \neq z_j \text{ if } i\neq j
	\end{array}
	\right\}.
\]
Let $J$ be an almost complex structure; we can then consider a bundle $\cE^{m-1,q}$ whose fiber above a point $u \in \cB^{m,q}_{A,0,k}$ is
\[
	\cE^{m-1,q}_{u} := W^{m-1,q} (S^2, \Lambda^{0,1}\otimes_J u^*TM).
\]
We view the Cauchy--Riemann operator $\delbarj := \tfrac{1}{2} (du + J \circ du \circ i)$ as a section of this bundle, i.e.,
\[
	\begin{tikzcd}
		\cE^{m-1,q} \arrow[r] & \cB^{m,q}_{A,0,k} \arrow[l, bend right, swap, "\delbarj"]
	\end{tikzcd}
\]
The linearization of this section is a Fredholm operator of index $2n+2c_1(A)+2k$.

We denote the set of $J$-holomorphic maps by
\[
	\widetilde{\cM}_{A,0,k}(J) :=
	\left\{
	\begin{array}{c}
		u: S^2 \to M \\
		\{z_1,\ldots,z_k\} \in S^2
	\end{array}
	\middle|
	\begin{array}{c}
		\tfrac{1}{2} (du + J \circ du \circ i) = 0 \\
		u_* [S^2] = A, z_i \neq z_j \text{ if } i \neq j
	\end{array}
	\right\}.
\]
It is identifiable with the zero set of the Cauchy--Riemann section, i.e., $\widetilde{\cM}_{A,0,k}(J) = \delbarjinv (0)$. It carries a natural reparametrization action by the group $\text{PSL}(2;\C)$. Consistent with our prior notation, we denote the quotient by
\[
	\cM_{A,0,k}(J) = \widetilde{\cM}_{A,0,k}(J) / \text{PSL}(2;\C).
\]
A $J$-holomorphic map is called \textbf{simple} if it is not a multiple cover. We denote the set of simple $J$-holomorphic maps by $\widetilde{\cM}^*_{A,0,k}(J)$ and we denote the set of simple $J$-holomorphic curves by $\cM^*_{A,0,k}(J) := \widetilde{\cM}^*_{A,0,k}(J) / \text{PSL}(2;\C)$.

Deep insight into this setup by Gromov in \cite{G} made possible the following:
\begin{itemize}
	\item The existence of regular almost complex structures such that the Cauchy--Riemann section is transverse to the zero section at all simple maps $u\in \cB^{m,q}_{A,0,k}$. Hence by the implicit function theorem the set $\widetilde{\cM}^*_{A,0,k}(J)$ is a manifold of dimension $2n+2c_1(A)+2k$. After taking the quotient by the reparametrization action the set $\cM^*_{A,0,k}(J) := \widetilde{\cM}^*_{A,0,k}(J) / \text{PSL}(2;\C)$ is a manifold of dimension $2n+2c_1(A) +2k -6$.
	\item The set $\cM_{A,0,k}(J)$ may be compactified in a natural way into $\CM_{A,0, k}(J)$ by adding additional nodal curves; this is the stable map compactification of Kontsevich.
\end{itemize}

The GW-moduli space $\CM_{A,0,k}(J)$ is therefore a compact topological space, while the subset of simple non-noded curves $\cM^*_{A,0,k}(J) \subset \CM_{A,0,k}(J)$ has the structure of a (usually noncompact) manifold.

\subsubsection{Pseudocycles}

Given a smooth manifold $M$, we say an arbitrary subset $B\subset M$ has \textbf{dimension at most} $d$, written as $\dim B \leq d$, if there exists a $d$-dimensional manifold $N$ and a smooth map $g: N \to M$ such that $B \subset g(N)$.

Given a continuous map between topological spaces, $h: A\to B$, we define the \textbf{omega limit set} of $h$ as follows:
\[	\Omega_h := \bigcap_{ \substack{K \subset A, \\ K \text{ compact}}}	\overline{h(A\setminus K)}.	\]
A point $x$ belongs to $\Omega_h$ if and only if $x$ is a limit point for a sequence $h(z_i)$ where $z_i \in A$ does not converge \cite{MSbook}*{Exer.~6.5.2. (i)}.
It follows that if $h$ can be continuously extended from $A$ to a larger compact topological space $\cA$ which contains $A$, then $\Omega_h \subset h(\cA \setminus A)$.

\begin{definition}
	Let $M$ be a smooth manifold. A smooth map
	\[f: V \to M\]
	is a $d$-dimensional \textbf{pseudocycle} if $V$ is an oriented $d$-dimensional manifold $V$ such that $f(V)$ has compact closure and such that $\dim \Omega_f \leq d - 2$.
\end{definition}

A pseudocycle $f:V \to M$ and a submanifold $X \subset M$ are \textbf{transverse}, written $f \pitchfork X$, if
\[
	\Omega_f\cap X = \emptyset
\]
and for all $x \in f^{-1}(X) \subset V$
\[
	Df_x (T_xV) \oplus T_{f(x)} X = T_{f(x)} M.
\]

Observe that if $f \pitchfork X$ then $f^{-1}(X)$ is a compact set, and hence when $\dim V + \dim X = \dim M$ then $f^{-1}(X)$ is also a finite set.
Thus we have a well-defined \textbf{intersection number}
\[
	f (V) \cdot X := \sum_{x \in f^{-1}(X)} \sign(x)
\]
where $\sign(x) = \pm 1$ is positive if $Df_x (T_xV) \oplus T_{f(x)} $ has the same orientation as $T_{f(x)} M$ and negative otherwise.
As it turns out, this intersection number depends only on the homology class of the submanifold and on the cobordism class of the pseudocycle.
Due to an oft cited result by Thom \cite{thom1954quelques}*{Thm.~II.1}, given a homology class in $H_*(M;\Q)$ we can find a submanifold which represents this class; we can therefore find a basis for $H_*(M;\Q)$ consisting of representing submanifolds $X_i$, and can therefore write an arbitrary homology class as $\sum_i k_i [X_i]$ for rational coefficients $k_i\in \Q$.

\begin{definition}[{\cite{MSbook}*{\S6.5}}]
	Given a pseudocycle $f: V\to M$ we can define a homomorphism by evaluating the intersection number with a basis of representing submanifolds and linear extension:
	\[
		\Phi_f : H_*(M;\Q) \to \Q, \qquad \sum_i k_i [X_i] \mapsto  \sum_i k_i (f(V) \cdot X_i).
	\]
\end{definition}

\subsubsection{The evaluation map on the Gromov--Witten moduli space is a pseudocycle}

The semipositive condition ensures the following: a nodal curve will not contain negative index domain components. Thus, we can control the dimension of the strata of the nodal tree compactification.

In the following situation the GW-moduli space has the structure of a pseudocycle.

\begin{theorem}[{\cite{MSbook}*{Thm.~6.6.1}}]
	Let $A\in H_2(M;\Z)$ such that
	\begin{equation}\label{eq:A_equals_mB}
		A=mB,\ c_1(B)=0 \quad \implies \quad m=1,
	\end{equation}
	for any $m\in \Z^{\geq 0}$ and $B$ a spherical homology class. Then there exists a regular almost complex structure $J$ such that the evaluation map
	\[
		ev : \cM^*_{A,0,k} (J) \to M^k
	\]
	is a pseudocycle of dimension $2n+2c_1(A)+2k -6$. Furthermore, for different choices of $J$ the resultant pseudocycles are cobordant.
\end{theorem}

In the general case, there is still the possibility of both non-noded multiple covers in the top stratum, as well as multiple covers which occur as domain components of a nodal curve in the lower strata. Moreover, we will also wish to incorporate homology classes from the GK-space. To define pseudocycle GW-invariants in full generality, additional constructions are necessary, we discuss this in \S\ref{sec:pseudocycle-full-gw-invariant}.

Assuming that \eqref{eq:A_equals_mB} holds, we may define the \textbf{pseudocycle GW-invariant} as the homomorphism
\begin{equation}\label{eq:pseudocycle-gw-invariant}
	\text{pseudocycle-}\GW_{A,0,k} : H_* (M;\Q)^{\otimes k} \to \Q
\end{equation}
via the intersection number of the pseudocycle $ev : \cM^*_{A,0,k} (J) \to M^k$ with a basis of representing submanifolds $\cX \subset M$:
\[
	\text{pseudocycle-}\GW_{A,0,k} ([\cX_1],\ldots,[\cX_k]) :=
	ev({\cM^*_{A,0,k}(J)}) \cdot \left(\cX_1 \times\cdots\times \cX_k \right).
\]
The invariant does not depend on the choice of regular almost complex structure $J$, nor on the choice of representing basis.


\subsection{Gromov--Witten invariants via polyfold theory}

We now recall the abstract perturbation scheme of polyfold theory used to define the polyfold GW-invariants.
Fix a homology class $A\in H_2(M;\Z)$, fix the genus $g=0$, and let $k\geq 0$ be the number of marked points.
Associated to $A,g,k$ we may write the \textbf{GW-polyfold} as the set
\begin{gather*}
	\cZ_{A,0,k} :=
	\left\{
	\begin{array}{c}
		u: S^2 \to M \\
		\{z_1,\ldots,z_k\} \in S^2
	\end{array}
	\middle|
	\begin{array}{c}
		u \in W^{2,3}(S^2,M), \\ u_* [S^2] = A,\\
		z_i \neq z_j \text{ if } i \neq j
	\end{array}
	\right\}
	\biggm/
	\begin{array}{c}
		\phi \in \text{PSL}(2,\C) \\
		\phi(z_i) = z'_i
	\end{array}
	\\
	\bigsqcup \ \{\text{stable nodal curves}\}.
\end{gather*}
The top stratum consists of all stable non-noded curves, as opposed to only the curves which are $J$-holomorphic.
It carries a natural second countable, paracompact, Hausdorff topology and a natural ``scale smooth'' structure.

We can define a strong polyfold bundle $\cW_{A,0,k}$ whose fiber over a stable curve consists of complex anti-linear forms. Over the top stratum of non-noded curves we can write simply
\[
	(\cW_{A,0,k})_{[(u,z_1,\ldots, z_k)]} : = W^{2,2}(S^2, \Lambda^{0,1}\otimes_J u^*TM) / \text{PSL}(2,\C).
\]
In this context we may also view the Cauchy--Riemann operator as a sc-smooth proper Fredholm section of this bundle,
\[
	\begin{tikzcd}
		\cW_{A,0,k} \arrow[r] & \cZ_{A,0,k}. \arrow[l, bend right, swap, "\delbarj"]
	\end{tikzcd}
\]
The GW-moduli space, already the quotient of groups of automorphisms and compactified with nodal $J$-holomorphic curves, then appears as the compact zero set of this section,
\[
	\CM_{A,0,k}(J) = \delbarjinv (0) \subset \cZ_{A,0,k}.
\]
\begin{theorem}[Polyfold regularization theorem, {\cite{HWZbook}*{Thm.~15.4, Cor.~15.1}}]
	There exists a class of $\ssc^+$-multisections $p$ such that the perturbed multisection $\delbarj+p$ is transverse to the zero section, and thus the perturbed zero set
	\[\CM_{A,0,k}(p):=(\delbarj + p)^{-1}(0)\]
	has the structure of a compact oriented ``weighted branched orbifold.''

	Furthermore, different choices of $\ssc^+$-multisections result in perturbed zero sets which are cobordant:
	given $\ssc^+$-multisections $p_1$, $p_2$ there exists a compact oriented weighted branched orbifold $\cB$ with boundary:
	\[
		\partial \cB = -\CM_{A,0,k}(p_1) \sqcup \CM_{A,0,k}(p_2).
	\]
\end{theorem}

\subsubsection{Weighted branched orbifolds at the local level}

We don't really want to give the complete definition of a weighted branched orbifold; to describe the global smooth structure we would need to introduce the modern language of ep-groupoids, and this would take us too far afield, see the exposition in \cite{HWZbook} or \cite{schmaltz2019steenrod} for details as needed.
In order to define intersection numbers it is enough to understand weighted branched orbifolds on a local level.

A \textbf{weighted branched orbifold} is a topological space $\cS$ which can be described locally as follows.
Given a point $x\in \cS$ there exists:
\begin{itemize}
	\item a collection of open subsets $U_i \subset \R^n$ called \textbf{local branches} together with positive rational weights $w_i$, indexed by a finite set $i\in I$,
	\item an ambient open set\footnote{This set has the structure of a $\ssc$-retract.} $V$ such that the local branches are properly embedded via inclusion maps
	      $U_i \hookrightarrow V;$
	      we identify $U_i$ with its image in $V$ and so we can write without ambiguity $U_i \subset V$ and $\cup_{i\in I} U_i \subset V$,
	\item a finite group $G$ which acts on $V$ by homeomorphisms and such that $\cup_{i\in I} U_i$ is invariant under this action,
	\item a quotient map $\pi: \cup_{i\in I} U_i  \to (\cup_{i\in I} U_i) / G$ defined as usual via the group action by $G$ on $\cup_{i\in I} U_i$,
	\item a local homeomorphism
	      \[
		      \phi: (\cup_{i\in I} U_i) / G \to \cS
	      \]
	      such that $\hat{x} := (\phi\circ \pi )^{-1}(x)$ lies on every branch, i.e., $\hat{x} \in U_i$ for all $i\in I$.
\end{itemize}
Let $M$ be a manifold. A continuous map $f:\cS \to M$ is smooth if and only if the composition on every local branch $f\circ \phi \circ \pi: U_i \to M$ is smooth.

We say that $f$ is \textbf{transverse} to a submanifold $X\subset M$, written $f \pitchfork X$, if for every $x \in f^{-1} (X) \subset \cS$, and for every $\hat{x}$ on some branch $U_i$ which maps via $\phi\circ \pi$ to $x$, we require
\[D(f\circ \phi \circ \pi)_{\hat{x}} U_i \oplus T_{f(x)} X = T_{f(x)} M.\]
By \cite{schmaltz2019steenrod}*{Lem.~3.12} when $\dim \cS + \dim X = \dim M$ the points $x \in f^{-1}(X)$ are isolated, and hence if $\cS$ is compact, $f^{-1}(X)$ consists of a finite set of points.

We may then define the \textbf{intersection number} as
\[
	f(\cS) \cdot X := \sum_{x \in f^{-1}(X)} \left(	\frac{1}{\sharp \bm{G}^{\textrm{eff}}(\hat{x})}\sum_{i\in I}	\sign (\hat{x}) w_i \right)
\]
where $\bm{G}^\text{non-eff} := \{	g \in \bm{G}	\mid	g (V) = V	\}$ and $\bm{G}^\text{eff} := \bm{G} / \bm{G}^\text{non-eff}.$
(It is of course possible that at a point $x \in \cS$ there is no branching and only a trivial group action; the intersection number at such a point reduces to the usual formula for an intersection number.)
As one would expect, this intersection number is independent of the cobordism class of the weighted branched orbifold.

As in the pseudocycle approach, this intersection number can be used to define a homomorphism on the rational homology $H_*(M; \Q)$ of a manifold $M$.
\begin{definition}[{\cite{schmaltz2019steenrod}*{Def.~3.13}}]
	Given a smooth map on a compact oriented weighted branched orbifold $f:\cS \to M$ we can define a homomorphism via intersection numbers with submanifolds as
	\[
		\Psi_f : H_* (M;\Q) \to \Q, \qquad \sum_i k_i [X_i] \mapsto \sum_i k_i (f(\cS) \cdot X_i).
	\]
\end{definition}

\subsubsection{The polyfold GW-invariants as intersection numbers}

Although originally defined in \cite{HWZGW} in terms of the branched integration of $\ssc$-smooth differential forms, in \cite{schmaltz2019steenrod} it was proved that the polyfold GW-invariants may equivalently be defined in terms of intersection numbers.

The \textbf{polyfold GW-invariant} is the homomorphism
\begin{equation}\label{eq:polyfold-gw-invariant}
	\text{polyfold-}\GW_{A,0,k} : H_* (M;\Q)^{\otimes k} \to \Q
\end{equation}
defined via the intersection number of a weighted branched orbifold with a basis of representing submanifolds $\cX \subset M$:
\[
	\text{polyfold-}\GW_{A,0,k} ([\cX_1],\ldots,[\cX_k]) :=
	ev({\CM_{A,0,k}(p)}) \cdot \left(\cX_1 \times\cdots\times \cX_k \right).
\]
The invariant does not depend on the choice of abstract perturbation, nor on the choice of representing basis.

We will also wish to consider GW-invariants with contributions from the homology groups of the GK-space. In the polyfold case, we may do this immediately as follows. We define the \textbf{full polyfold GW-invariant} as the homomorphism
\begin{equation}\label{eq:full-polyfold-gw-invariant}
	\text{polyfold-}\GW_{A,0,k} : H_* (M;\Q)^{\otimes k} \otimes H_* (\dmspace_{0,k}; \Q) \to \Q
\end{equation}
via the intersection number of a weighted branched orbifold with a basis of representing submanifolds $\cX \subset M$ and with a basis of representing submanifolds $\cB \subset \dmspace_{0,k}$:
\[
	\text{polyfold-}\GW_{A,0,k} ([\cX_1],\ldots,[\cX_k];[\cB]) :=
	(ev\times \pi )({\CM_{A,0,k}(p)}) \cdot \left(\cX_1 \times\cdots\times \cX_k \times \cB\right).
\]
The invariant does not depend on the choice of abstract perturbation, nor on the choice of representing bases. It agrees with \eqref{eq:polyfold-gw-invariant} in the case that $\cB = \dmspace_{0,k}$.


\section{The abstract comparison theorem}
\label{sec:the_abstract_comparison_thm}

\subsection{Pseudocycles versus weighted branched orbifolds}

A pseudocycle has a top stratum covered by a single manifold; the dimensions of the lower strata are controllable, but there is, \textit{a priori}, no smooth structure in neighborhoods of these strata.

In contrast, weighted branched orbifolds may contain possible branching on the top stratum; the lower strata may contain orbifold singularities in addition to branching phenomena. Crucially, weighted branched orbifolds are smooth geometric spaces and the lower strata have local smooth structure everywhere; it is possible to consistently define notions of integration as in \cite{HWZint} or intersection number as in \cite{schmaltz2019steenrod}.

Pseudocycles and weighted branched orbifolds are fundamentally different geometric objects.
Equality of the pseudocycle and polyfold GW-invariants should not be taken for granted.

\begin{figure}[ht]
	\centering
	\includegraphics[width=\textwidth]{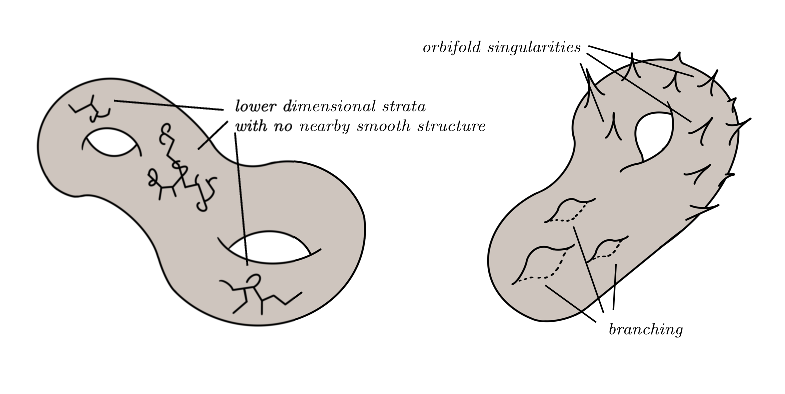}
	\caption{A pseudocycle versus a weighted branched orbifold}
\end{figure}

\subsection{Statement and proof of the abstract comparison theorem}

The pseudocycle and polyfold approach to regularization mesh together in the following description.
The Cauchy--Riemann operator defines a $\ssc$-smooth proper Fredholm section of the GW-polyfold:
\[
	\begin{tikzcd}
		\cW_{A,0,k} \arrow[r] & \cZ_{A,0,k}. \arrow[l, bend right, swap, "\delbarj"]
	\end{tikzcd}
\]
In the case $A=mB,\ c_1(B)=0 \implies m=1$, there exists an almost complex structure $J$ such that $\delbarj$ is transverse to the zero section for all simple non-noded stable curves; hence $\cM^*_{A,0,k}(J)$ is a manifold of dimension $\ind \delbarj = 2n+2c_1(A)+2k -6$. Furthermore, the image of the evaluation map restricted to the nonsimple/nodal stable curves has dimension at most $\ind \delbarj - 2$.
Abstracted, we offer the following theorem.

\begin{theorem}[Abstract comparison theorem]\label{thm:comparison-thm}
	Consider a strong polyfold bundle $\cW$ over a polyfold $\cZ$ and a $\ssc$-smooth proper Fredholm section $s$:
	\[
		\begin{tikzcd}
			\cW \arrow[r] & \cZ. \arrow[l, bend right, swap, "s"]
		\end{tikzcd}
	\]
	Denote the compact solution set as $\cM := s^{-1}(0)$.

	Let $M$ be a smooth manifold and let $f: \cZ \to M$ be a $\ssc$-smooth map.
	Suppose that there exists an open subset $\cZ^* \subset \cZ$ such that:
	\begin{itemize}
		\item $s|_{\cZ^*}\pitchfork 0$,
		\item $\dim f(\cM \setminus \cM^*) \leq \ind s -2$, where $\cM^* := \cM \cap \cZ^*$.
	\end{itemize}
	(This directly implies that there exists an open subset $\cZ^* \subset \cZ$ such that $f: \cM^* \to M$, where $\cM^* := \cM \cap \cZ^*$, is a pseudocycle of dimension $\ind s$.)

	Then, given any submanifold $X \subset M$ of dimension $\dim X = \dim M - \dim \cM^*$ there exists an abstract perturbation $p$ of the section $s$ such that the intersection number of the pseudocycle $f:\cM^* \to M$ with $X$ and the intersection number of the compact weighted branched orbifold $\cM(p)$ with $X$ are equal, i.e.,
	\[
		\underbrace{f(\cM^*)}_{\text{pseudocycle}} \cdot\ X
		\quad=
		\underset{\substack{  \text{compact weighted}\quad \\ \text{branched orbifold}\quad }}{f(\underbrace{\cM(p)}) \cdot X.}
	\]
\end{theorem}
\begin{proof}
	Without loss of generality, we may assume that we have perturbed $X$ so that it is transverse to $\cM^*$ and such that $X\cap f(\cM \setminus \cM^*) = \emptyset$.

	Note that $f^{-1}(X)$ and $\cM \setminus \cM^*$ are closed subsets of $\cZ$; and moreover that $f^{-1} (X)\cap (\cM \setminus \cM^*) = \emptyset$. Since polyfolds are normal topological spaces, we can find a neighborhood $\cU$ of $\cM \setminus \cM^*$ such that $f^{-1}(X) \cap \cU = \emptyset$.

	By assumption, $s$ is already transverse on the compliment of $\cU$; hence we may restrict the support of any regularizing abstract perturbation to $\cU$.
	Any possible branching/isotropy phenomena of the resulting weighted branched orbifold $\cM(p)$ will therefore be restricted to $\cU$, and moreover $\cM(p)$ and $\cM^*$ are identical outside of $\cU$. It follows that all points of intersection with $X$ must lie outside of $\cU$. Hence the points of intersection for $\cM^*$ with $X$ and $\cM(p)$ with $X$ are equal, and hence the intersection numbers for $\cM^*$ considered as a pseudocycle or $\cM(p)$ considered as a weighted branched orbifold are also equal.
\end{proof}

\begin{figure}[ht]
	\centering
	\begin{subfigure}[b]{0.3\textwidth}
		\centering
		\includegraphics[width=\textwidth]{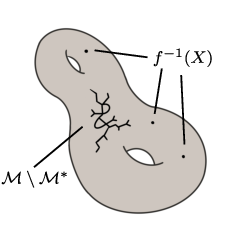}
		\caption{The preimage $f^{-1}(X)$ is disjoint from the lower dimensional strata $\mathcal{M} \setminus \mathcal{M}^*$.\\ \mbox{} \\}
	\end{subfigure}
	\hfill
	\begin{subfigure}[b]{0.3\textwidth}
		\centering
		\includegraphics[width=\textwidth]{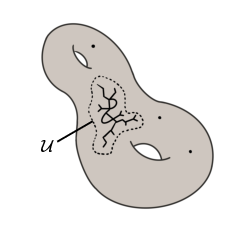}
		\caption{We may find an open neighborhood $\mathcal{U}$ of the lower dimensional strata, such that the neighborhood remains disjoint from $f^{-1}(X)$.}
	\end{subfigure}
	\hfill
	\begin{subfigure}[b]{0.3\textwidth}
		\centering
		\includegraphics[width=\textwidth]{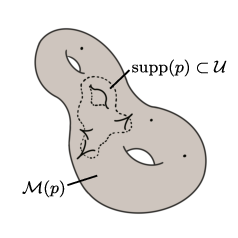}
		\caption{We may choose a regular perturbation $p$ with support contained within $\mathcal{U}$. The points of intersection of $\mathcal{M}(p)$ and $X$ are still $f^{-1}(X)$.}
	\end{subfigure}
	\caption{Picture proof of the abstract comparison theorem}
\end{figure}

We may immediately apply the abstract comparison theorem to see that for any collection of submanifolds $\cX_i \subset M$ such that $\sum_{i=1}^k (2n - \dim \cX_i) = 2n+2c_1(A)+2k -6 $ we have
\[
	\underbrace{ev(\cM^*_{A,0,k}(J))}_{\text{pseudocycle}} \cdot (\cX_1 \times \cdots \times \cX_k)
	\quad= \quad
	\underset{\substack{  \text{compact weighted} \\ \text{branched orbifold} }}{ev(\underbrace{\CM_{A,0,k}(p)})} \cdot (\cX_1 \times \cdots \times \cX_k)
\]
and thus we obtain the following corollary.

\begin{corollary}
	Let $(M,\ww)$ be a semipositive symplectic manifold.
	Let $A\in H_2(M;\Z)$ such that
	\[
		A=mB,\ c_1(B)=0 \quad \implies \quad m=1.
	\]
	for $m\in \Z^{\geq 0}$ and $B$ a spherical homology class. Then the pseudocycle \eqref{eq:pseudocycle-gw-invariant} and polyfold Gromov--Witten invariants \eqref{eq:polyfold-gw-invariant} for the class $A$ are equal,
	\[
		\text{pseudocycle-}\GW_{A,0,k} = \text{polyfold-}\GW_{A,0,k}.
	\]
\end{corollary}


\section{The pseudocycle approach to defining the full Gromov--Witten invariants}
\label{sec:pseudocycle-full-gw-invariant}

In order to completely define the full GW-invariant via the pseudocycle approach, we must now deal with two problems:

\begin{problem}\label{prob:problem1}
The case $A=mB$, $c_1(B)=0$, $m >1$. In this case, it is no longer possible to control the dimension of the nonsimple non-noded elements in the GW-moduli space.
\end{problem}

\begin{problem}\label{prob:problem2}
Contributions from homology classes $B\in H_*(\dmspace_{0,k};\Q)$. In general, it is not possible to show that the map
\[
	ev \times \pi : \CM_{A,0,k} (J) \to M^k \times \dmspace_{0,k}
\]
is a pseudocycle, so the natural approach of defining the GW-invariant as an intersection number with a submanifold $\cX_1 \times \cdots \times \cX_k \times \cB$ is not possible.
\end{problem}

The pseudocycle approach deals with these problems with the following solutions:

\begin{solution}
	Allow more general perturbation of the almost complex structure. In particular, allow the almost complex structure to vary with the domain, i.e., consider families of the form $\{J_z\}_{z\in S^2}$.
\end{solution}
\begin{solution}
	Restrict the GW-moduli space to fixed tree types $T$, corresponding to submanifolds of the GK-space $\cB_T\subset \dmspace_{0,k}$ which generate the homology of $H_*(\dmspace_{0,k};\Q)$. The restriction to such a tree yields a GW-moduli space $\CM_T(J)$ identifiable with the preimage of $\cB_T$ via $\pi: \CM_{A,0,k}(J) \to \dmspace_{0,k}$, i.e., $\CM_T(J) = \pi^{-1} (\cB_T)$.
\end{solution}

\subsection{Domain dependent perturbation and graphs}
\label{subsec:domain-dependent-perturbation}

Consider a smooth family $\{J_z\}_{z\in S^2}$ of $\ww$-compatible almost complex structures which depend on a point of the domain $S^2$. Solutions of the Cauchy--Riemann operator become
\begin{equation}\label{eq:domain-dependent-cr}
	\tfrac{1}{2}(du + J_z(u(z)) \circ du \circ i) = 0.
\end{equation}

It is meaningful to consider solutions to a domain dependent almost complex structure in a global setting as \textit{graphs}.
Define a product manifold and product almost complex structure as follows:
\[
	\widetilde{M} := S^2 \times M, \qquad \widetilde{J}(z,x) : = i \oplus J_z(x).
\]
We then observe that $u: S^2 \to M$ is a solution of \eqref{eq:domain-dependent-cr} if and only if its graph
\[\widetilde{u}: S^2 \to \widetilde{M},\qquad \widetilde{u}(z) :=(z,u(z))\]
is $\widetilde{J}$-holomorphic, i.e., $\delbar_{\widetilde{J}} \widetilde{u} =0$.

Let $\widetilde{A} := [S^2 \times \text{pt}] + \iota_* A \in H_2(\widetilde{M};\Z)$ where $A\in H_2(M;\Z)$ and $\iota : M \to \widetilde{M},\ x \mapsto (\text{pt},x)$. We may consider the associated GW-moduli space
$\CM_{\widetilde{A},0,k} (\widetilde{J}).$
The usual proofs of transversality and compactness carry forward, even with the additional requirement that one consider almost complex structures on $\widetilde{M}$ which split with respect to the factors, i.e., are of the form $\widetilde{J}:= i \oplus J_z$.
Notice however that the expected dimension of the graph GW-moduli space $\CM_{\widetilde{A},0,k} (\widetilde{J})$ is $2n+2c_1(A) +2k$ while the expected dimension of the usual GW-moduli space $\CM_{A,0,k}(J)$ is $2n+2c_1(A)+2k -6$.
To define GW-invariants via graphs in such a way that they will also agree with the pseudocycle GW-invariants \eqref{eq:pseudocycle-gw-invariant} it is necessary to fix the first three marked points; on the top stratum we may do this directly and write
\[
	\cM_{\widetilde{A},0,k}^{(0,1,\infty)} (\widetilde{J})
	:=
	\left\{
	\begin{array}{c}
		\tilde{u}: S^2 \to \widetilde{M} \\
		\{z_1,\ldots,z_k\} \in S^2
	\end{array}
	\middle|
	\begin{array}{c}
		\tfrac{1}{2} (d\tilde{u} + \widetilde{J} \circ d\tilde{u} \circ i) = 0 \\
		\tilde{u}_* [S^2] = \widetilde{A},                                     \\
		z_i \neq z_j \text{ if } i \neq j,                                     \\
		z_1 =0, \ z_2 = 1,\ z_3 = \infty
	\end{array}
	\right\}
	\biggm/
	\text{PSL}(2,\C).
\]
Consider the map
\begin{equation}\label{eq:theta-fixes-marked-points}
	\Theta : \CM_{\widetilde{A},0,k} (\widetilde{J})\xrightarrow{ev_1\times ev_2 \times ev_3} \widetilde{M}\times \widetilde{M}\times \widetilde{M} \xrightarrow{p_1\times p_1 \times p_1} S^2 \times S^2 \times S^2.
\end{equation}
where $p_1 : \widetilde{M} = S^2 \times M \to S^2$ is the projection onto the first factor.
The natural compactification of $\cM_{\widetilde{A},0,k}^{(0,1,\infty)} (\widetilde{J})$ is precisely equal to the preimage of $(0,1,\infty) \in S^2\times S^2\times S^2$, i.e.,
\[
	\CM_{\widetilde{A},0,k}^{(0,1,\infty)} (\widetilde{J}) = \Theta^{-1} (0,1,\infty).
\]
Finally, observe that for a generic almost complex structure $\widetilde{J}$ the top stratum of the GW-moduli spaces will consist only of simple curves; the domain dependent perturbation generically kills the symmetry for multiple covers. In other words, generically we have $\cM_{\widetilde{A},0,k}^* (\widetilde{J}) = \cM_{\widetilde{A},0,k} (\widetilde{J})$ and $\cM_{\widetilde{A},0,k}^{*(0,1,\infty)} (\widetilde{J}) = \cM_{\widetilde{A},0,k}^{(0,1,\infty)} (\widetilde{J})$.

\begin{theorem}{\cite{MSbook}*{Thm.~6.7.1}}
	There exists a regular almost complex structure on $\widetilde{M}$ of the form $\widetilde{J}= i \oplus J_z$ such that the evaluation map
	\[
		ev:\cM_{\widetilde{A},0,k}^{(0,1,\infty)} (\widetilde{J}) \to M^k
	\]
	is a pseudocycle of dimension $2n+2c_1(A)+2k-6$. Furthermore, for different choices of regular almost complex structure the resultant pseudocycles are cobordant.

	Finally, in the case \eqref{eq:A_equals_mB}, i.e., \[A =mB,\ c_1(B)=0 \implies m=1\] for any $m\in \Z^{\geq 0}$, the pseudocycle $\cM_{\widetilde{A},0,k}^{(0,1,\infty)} (\widetilde{J})$ is cobordant to the usual GW-pseudocycle $\cM^*_{A,0,k}(J)$.
\end{theorem}

It is important to understand both the significance and reasoning behind this last statement, i.e., that these differently defined GW-pseudocycles are cobordant.
Assuming the prior hypothesis \eqref{eq:A_equals_mB}, and assuming that $J$ is regular, we may define a diffeomorphism between the top stratum of the moduli spaces (which are finite-dimensional manifolds) as follows:
\begin{align}
	\Psi: \cM^*_{A,0,k}(J) & \to \cM_{\widetilde{A},0,k}^{*(0,1,\infty)} (i\oplus J) \label{eq:inclusion-map-top-stratum-gw-moduli} \\
	[(u,z_1,\ldots,z_k)]   & \mapsto [\left((\id, u\circ \phi^{-1} ), \phi(z_1),\ldots, \phi(z_k) \right)] \nonumber
\end{align}
where $\phi :S^2 \to S^2$ is the unique M\"obius transformation such that $\phi(z_1) = 0,\ \phi(z_2) = 1,\ \phi(z_3) = \infty$.
Moreover, this map $\Psi$ commutes with the evaluation maps, thus giving an identification between these moduli spaces considered as pseudocycles.
From here, one can then choose a path between $i \oplus J$ and any regular $i\oplus J_z$ to obtain a cobordism between pseudocycles, justifying the last statement above.

Using the above theorem, we may take the previously defined pseudocycle GW-invariant \eqref{eq:pseudocycle-gw-invariant} and extend it to all spherical homology classes $A\in H_2(M;\Z)$ via the graph GW-moduli spaces by taking the intersection number of the pseudocycle $ev:\cM_{\widetilde{A},0,k}^{(0,1,\infty)} (\widetilde{J}) \to M^k$ with a basis of representing submanifolds $\cX \subset M$:
\begin{equation}\label{eq:pseudo-graph-gw}
	\text{pseudocycle-}\GW_{A,0,k} ([\cX_1],\ldots,[\cX_k]) :=
	ev(\cM_{\widetilde{A},0,k}^{(0,1,\infty)} (\widetilde{J})) \cdot \left(\cX_1 \times\cdots\times \cX_k \right).
\end{equation}

\subsection{Fixing the tree type via a submanifold of the Grothendieck--Knudsen space}
\label{subsec:fixing_the_tree_type}

In order to allow for contributions to the GW-invariant from homology classes $B\in H_*(\dmspace_{0,k};\Q)$ we now define GW-invariants for fixed tree type.

Let $k\geq 3$. A \textbf{$k$-labeled tree} consists of a tree $T$ (which consists of a set of vertices $V$ and edges $E$) together with a labeling of the vertices, $\Lambda : \{1,\ldots, k\} \to V$.
A $k$-labeled tree is \textbf{stable} if each vertex $\alpha \in T$ satisfies
\[
	\# \Lambda^{-1}(\alpha) + \# \{\beta \in T \mid \alpha E \beta \} \geq 3,
\]
in other words, the number of labels on $\alpha$ and the number of edges with a vertex at $\alpha$ is at least $3$.

We can restrict the GK-space $\dmspace_{0,k}$ to the subset of stable noded Riemann surfaces which are modeled on such a stable $k$-labeled tree $T$; we denote this as $\cM_{0,T}$. The natural compactification, $\CM_{0,T}$, is a closed submanifold of $\dmspace_{0,k}$.
Keel showed that these submanifolds generate the rational homology of $\dmspace_{0,k}$.

\begin{theorem}[{\cite{keel1992intersection}}]\label{thm:homology-gk-space}
	The homology classes $[\CM_{0,T}]$ form a basis for the rational homology groups $H_*(\dmspace_{0,k};\Q)$.
\end{theorem}

We may restrict the GW-moduli space to fixed tree types as follows. Consider the projection map
\[
	\pi: \CM^{(0,1,\infty)}_{\widetilde{A},0,k}(\widetilde{J}) \to \dmspace_{0,k}
\]
and consider the preimage
\[
	\cM^{(0,1,\infty)}_{\widetilde{A},0,k} \big|_T (\widetilde{J}):= \pi^{-1}(\cM_{0,T}).
\]
One may check that the natural Gromov compactification of this moduli space, $\CM^{(0,1,\infty)}_{\widetilde{A},0,k} \big|_T (\widetilde{J})$, is compatible with compactification of $\CM_{0,T}$ via the projection map in the sense that $\CM^{(0,1,\infty)}_{\widetilde{A},0,k} \big|_T (\widetilde{J})= \pi^{-1}(\CM_{0,T})$.
We call $\CM^{(0,1,\infty)}_{\widetilde{A},0,k} \big|_T (\widetilde{J})$ the \textbf{GW-moduli space of stable curves modeled on a stable $k$-labeled tree $T$}.

\begin{theorem}{\cite{MSbook}*{Thm.~6.7.11}}\label{thm:tree-pseudocycle}
	There exists a regular almost complex structure of the form $\widetilde{J}=i\oplus J_z$ such that the evaluation map
	\[
		ev: \cM^{(0,1,\infty)}_{\widetilde{A},0,k} \big|_T (\widetilde{J}) \to M^k
	\]
	is a pseudocycle of dimension $2n+2c_1(A)+2k -6 - 2e(T)$, where $e(T)$ is the number of edges of the tree $T$. Furthermore, for different choices of $\widetilde{J}$ the resultant pseudocycles are cobordant.
\end{theorem}

Combining theorems~\ref{thm:homology-gk-space} and \ref{thm:tree-pseudocycle} we may define the \textbf{full pseudocycle GW-invariant} as the homomorphism
\begin{equation}\label{eq:full-pseudo-gw-invariant}
	\text{pseudocycle-}\GW_{A,0,k} : H_* (M;\Q)^{\otimes k} \otimes H_* (\dmspace_{0,k}; \Q) \to \Q
\end{equation}
by considering a basis of representing submanifolds $\cB_T \subset \dmspace_{0,k}$ of fixed tree type, and then evaluating the intersection number of the pseudocycle $ev : \cM^{(0,1,\infty)}_{\widetilde{A},0,k} \big|_T (\widetilde{J}) \to M^k$ with a basis of representing submanifolds $\cX \subset M$:
\[
	\text{pseudocycle-}\GW_{A,0,k} ([\cX_1],\ldots,[\cX_k]; [\cB_T]) :=
	ev\left(\cM^{(0,1,\infty)}_{\widetilde{A},0,k} \big|_T (\widetilde{J})\right) \cdot \left(\cX_1 \times\cdots\times \cX_k \right).
\]
The invariant does not depend on the choice of regular almost complex structure $\widetilde{J}$.

The geometric meaning of this definition is essentially the same as before; the GW-invariant is the signed count of the finite points of intersection of the GW-moduli space with the submanifolds $\cX_i \subset M$ and the submanifold $\cB_T\subset \dmspace_{0,k}$.
Indeed, observe the finite sets of points of intersection are identical, i.e.,
\begin{align}
	 & (ev\times\pi)\left(\cM^{(0,1,\infty)}_{\widetilde{A},0,k} (\widetilde{J})\right) \cap \left(\cX_1 \times\cdots\times \cX_k\times \cB_T \right) \label{eq:not_a_pseudocycle} \\
	 & \qquad = ev\left(\cM^{(0,1,\infty)}_{\widetilde{A},0,k} \big|_T (\widetilde{J})\right) \cap \left(\cX_1 \times\cdots\times \cX_k \right) \nonumber
\end{align}
which is due to the identification
\[
	\cM^{(0,1,\infty)}_{\widetilde{A},0,k} \big|_T (\widetilde{J}) = \pi^{-1}\left(  \cM^{(0,1,\infty)}_{\widetilde{A},0,k} (\widetilde{J}) \right).
\]
However as we observed in Problem~\ref{prob:problem2}, we cannot expect that $ev \times \pi$ is a pseudocycle and hence the intersection number for the left-hand side of equation~\eqref{eq:not_a_pseudocycle} does not make sense. Restricting our moduli spaces to trees allows us to consider a situation where a pseudocycle intersection number can be defined.

\section{The Gromov--Witten polyfold of graphs of stable curves}
\label{sec:polyfold_of_graphs}

In the previous section we examined the pseudocycle approach to defining the full GW-invariants.
In this section we mimic these constructions in the context of polyfold theory.

\subsection{The Gromov--Witten polyfold of graphs}

We now recast the global setup of the ``Cauchy--Riemann section as a graph'' in the language of polyfolds.
As before, we consider a product symplectic manifold with a spherical homology class
\[
	\widetilde{M} = S^2 \times M, \qquad \widetilde{A} = [S^2 \times \text{pt}] + \iota_* A \in H_2(\widetilde{M};\Z).
\]
The GW-polyfold of graphs is simply the usual GW-polyfold associated to this data, i.e.,
$\cZ_{\widetilde{A},0,k}$. For a product almost complex structure $\widetilde{J} = i \oplus J_z$ we can also define a strong polyfold bundle and use the Cauchy--Riemann operator to define a $\ssc$-smooth proper Fredholm section $\delbar_{\widetilde{J}}$ of index
\[
	\ind \delbar_{\widetilde{J}} = 2n +2c_1(A) + 2k.
\]
As in \S\ref{subsec:domain-dependent-perturbation}, we will fix the first three marked points in order to get comparable GW-invariants.
To do this, we extend the map $\Theta$ from equation~\eqref{eq:theta-fixes-marked-points} to the GW-polyfold, and so consider the composition
\[
	\hat{\Theta}: \cZ_{\widetilde{A},0,k} \xrightarrow{ev_1\times ev_2 \times ev_3} \widetilde{M}\times \widetilde{M}\times \widetilde{M} \xrightarrow{p_1\times p_1\times p_1} S^2 \times S^2 \times S^2.
\]
By \cite{schmaltz2019axioms}*{Prop.~3.3}, evaluation maps on GW-moduli spaces are $\ssc$-smooth submersions; the projection to the first factor $p_1 : \widetilde{M}= S^2\times M \to S^2$ is also a smooth submersion. Therefore the map $\hat{\Theta}$ is a $\ssc$-smooth submersion.
In order to give the preimage
\begin{equation}\label{eq:GW-polyfold_with_fixed_marked_points}
	\cZ_{\widetilde{A},0,k}^{(0,1,\infty)}: = \hat{\Theta}^{-1}(0,1,\infty)
\end{equation}
a smooth structure we may invoke a result of Filippenko, who proved a general transverse preimage theorem for polyfold theory.
Thus by \cite{filippenko2018constrained}*{Thm.~1.5}, without loss of generality\footnote{Technically, to get a polyfold structure on a transverse preimage we must concern ourselves with the grading of polyfolds and must shift the levels of the polyfold up by 1. While important to be aware of such modifications to the abstract framework this does not alter our current discussion.}, we obtain a polyfold structure on $\cZ_{\widetilde{A},0,k}^{(0,1,\infty)}$, which we call the \textbf{GW-polyfold of graphs with three fixed marked points}.
By this theorem, we also obtain a strong polyfold bundle over $\cZ_{\widetilde{A},0,k}^{(0,1,\infty)}$ (given by the restriction to this subpolyfold) and a $\ssc$-smooth proper Fredholm section $\delbar_{\widetilde{J}}|_{\cZ_{\widetilde{A},0,k}^{(0,1,\infty)}}$ of index
\[
	\ind \delbar_{\widetilde{J}}|_{\cZ_{\widetilde{A},0,k}^{(0,1,\infty)}} = 2n +2c_1(A)+2k-6.
\]
After perturbation of this section, we obtain a \textbf{perturbed graph GW-moduli space} $\CM_{\widetilde{A},0,k}^{(0,1,\infty)}(p)$ which has the structure of a compact oriented weighted branched orbifold.
We may then define an analogous \textbf{polyfold graph GW-invariant} via the intersection number with this perturbed graph GW-moduli space:
\begin{equation}\label{eq:polyfold-graph-gw-invariants}
	\text{polyfold-graph-}\GW_{A,0,k} ([\cX_1],\ldots,[\cX_k]) :=
	ev({\CM_{\widetilde{A},0,k}^{(0,1,\infty)}(p)}) \cdot \left(\cX_1 \times\cdots\times \cX_k \right);
\end{equation}
we may allow contributions from homology classes in $\dmspace_{0,k}$ and thus consider also:
\begin{align}\label{eq:full-polyfold-graph-gw-invariants}
	 & \text{polyfold-graph-}\GW_{A,0,k} ([\cX_1],\ldots,[\cX_k]; [\cB])                                                                                         \\
	 & \qquad	:=	(ev\times \pi)\left({\CM_{\widetilde{A},0,k}^{(0,1,\infty)}(p)}\right) \cdot \left(\cX_1 \times\cdots\times \cX_k \times \cB \right). \nonumber
\end{align}
We may immediately apply the abstract comparison theorem~\ref{thm:comparison-thm} to see that
\[
	\underbrace{ev(\cM^{*(0,1,\infty)}_{\widetilde{A},0,k}(\widetilde{J}))}_{\text{pseudocycle}} \cdot (\cX_1 \times \cdots \times \cX_k)
	\quad= \quad
	\underset{\substack{  \text{compact weighted} \\ \text{branched orbifold} }}{ev\big(\underbrace{\CM^{(0,1,\infty)}_{\widetilde{A},0,k}(p)}\big)} \cdot (\cX_1 \times \cdots \times \cX_k)
\]
and thus obtain the following corollary.

\begin{corollary}
	For any homology class $A\in H_2(M;\Z)$ the extended pseudocycle GW-invariants \eqref{eq:pseudo-graph-gw} and the polyfold graph GW-invariants \eqref{eq:polyfold-graph-gw-invariants} are equal.
\end{corollary}

\subsection{Fixing the tree type of the Gromov--Witten polyfold of graphs}

We now discuss how to fix the tree type of the graph GW-polyfold, mimicking the setup of \S\ref{subsec:fixing_the_tree_type}.
Given a stable $k$-labeled tree $T$, we associate the submanifold $\CM_{0,T} \subset \dmspace_{0,k}$, and consider the projection map
\[
	\pi : \cZ_{\widetilde{A},0,k}^{(0,1,\infty)} \to \dmspace_{0,k};
\]
we may restrict the GW-polyfold to the tree $T$ by taking the preimage of $\CM_{0,T}$:
\begin{equation}\label{eq:GW-polyfold_restricted_to_tree}
	\cZ_{\widetilde{A},0,k}^{(0,1,\infty)} \big|_T := \pi^{-1} (\CM_{0,T}).
\end{equation}
Once we give this a polyfold structure, we will call this the \textbf{GW-polyfold of graphs with three fixed marked points and which is restricted to a fixed tree}.

Ideally we would give this subset a polyfold structure via the transverse preimage theorem of Filippenko as above, however this is not possible because in general the projection map is not a submersion, nor is it transverse to the submanifold $\CM_{0,T}$, see \cite{schmaltz2019axioms}*{Prob.~4}. The linearization of $\pi$ vanishes on the tangent planes of the gluing parameters, while the tangent planes to the gluing parameters are also complementary to the tangent spaces of $\CM_{0,T}$. Hence we cannot in general obtain transversality without perturbing the submanifold $\CM_{0,T}$.

We may still give $\cZ_{\widetilde{A},0,k}^{(0,1,\infty)} \big|_T$ a polyfold structure, however.
Due to the explicit structure of the submanifold $\CM_{0,T}$, one can similarly construct a polyfold structure on $\cZ_{\widetilde{A},0,k}^{(0,1,\infty)} \big|_T$. For example, local uniformizers may be constructed explicitly by requiring that the gluing parameters which correspond to the edges of the tree $T$ are $0$.

Another possible approach to give $\cZ_{\widetilde{A},0,k}^{(0,1,\infty)} \big|_T$ a polyfold structure is the method used to prove \cite{schmaltz2019axioms}*{\S2.3, Thm.~4.7}.
The idea is that we consider a disjoint union of polyfolds modeled on fixed trees but without any incidence relation between the vertices of an edge. We then take the preimage of the diagonal under the evaluation map, which then enforces the desired incidence relation. Moreover, this preimage will have the structure of a polyfold by the transverse preimage theorem.
Explicitly, we consider a $k$-labeled tree $\widetilde{T}$ and a decomposition of $\widetilde{A} = \sum_{\alpha \in \widetilde{T}} A_\alpha$ such that each vertex $\alpha \in T$ satisfies the usual GW-stability conditions:
\[
	\# (\text{marked points, edges at } \alpha) \geq 3 \quad \text{or} \quad \langle[\ww],A_\alpha\rangle >0.
\]
We require that after forgetting the components $A_\alpha$ and contracting unstable vertices, we obtain our original tree $T$.
To each vertex we then associate a GW-polyfold $\cZ_{A_\alpha,0,k_\alpha}$ where $k_\alpha$ is the number of labels at the vertex $\alpha$ in addition to the number of edges at the vertex $\alpha$.
We then consider the map
\[
	\prod_{\alpha E \beta} (ev_{k_\alpha} \times ev_{k_\beta}) : \bigsqcup_{\ast \in \widetilde{T}} \cZ_{A_\ast,0,k_\ast} \to  (M\times M)^{\# E}
\]
We can recover the incidence relation for the edges by taking the preimage of the diagonal $\Delta^{\#E} \subset (M\times M)^{\# E}$.
Thus,
\[
	\cZ_{\widetilde{A},0,k}^{(0,1,\infty)} \big|_T = \left(\prod_{\alpha E \beta} (ev_{k_\alpha} \times ev_{k_\beta}) \right)^{-1} \left(\Delta^{\#E}\right)
\]
As we have noted, evaluation maps are $\ssc$-smooth submersions, and hence by the transverse preimage theorem we conclude that $\cZ_{\widetilde{A},0,k}^{(0,1,\infty)} \big|_T$ has a polyfold structure.

Either way, $\cZ_{\widetilde{A},0,k}^{(0,1,\infty)} \big|_T$ is a polyfold, the Cauchy--Riemann operator restricted to $\cZ_{\widetilde{A},0,k}^{(0,1,\infty)} \big|_T$ is a $\ssc$-smooth proper Fredholm section of index $2n+2c_1(A)+2k - 6 -e(T)$, and we may define polyfold graph GW-invariants for the fixed tree $T$ via the intersection number:
\begin{align}\label{eq:polyfold-graph-tree-gw-invariants}
	 & \text{polyfold-fixed-tree-}\GW_{A,0,k} ([\cX_1],\ldots,[\cX_k];[\CM_{0,T}])                                                                \\
	 & \qquad	:=	ev\left({\CM_{\widetilde{A},0,k}^{(0,1,\infty)}\big|_T (p)}\right) \cdot \left(\cX_1 \times\cdots\times \cX_k \right). \nonumber
\end{align}
Once again applying the abstract comparison theorem~\ref{thm:comparison-thm} we see that
\[
	\underbrace{ev(\cM^{(0,1,\infty)}_{\widetilde{A},0,k}\big|_T (\widetilde{J}))}_{\text{pseudocycle}} \cdot (\cX_1 \times \cdots \times \cX_k)
	\quad= \quad
	\underset{\substack{  \text{compact weighted} \\ \text{branched orbifold} }}{ev \big(\underbrace{{\CM_{\widetilde{A},0,k}^{(0,1,\infty)}\big|_T (p)}}\big)} \cdot (\cX_1 \times \cdots \times \cX_k)
\]
and thus obtain the following corollary.

\begin{corollary}
	The full pseudocycle GW-invariant \eqref{eq:full-pseudo-gw-invariant} is equal to the polyfold GW-invariant defined via fixed tree types \eqref{eq:polyfold-graph-tree-gw-invariants}.
\end{corollary}


\section{Equality of the various polyfold Gromov--Witten invariants}
\label{sec:equality_polyfold_GW_invariants}

At this point, we have entered a pattern.
In order to define a GW-invariant in all situations and with contributions from the homology classes of the GK-space, the pseudocycle approach has made tweaks and modifications to the definition of a GW-invariant.
At every step, we have mimicked these modified constructions in the language of polyfold theory, and used the abstract comparison theorem to prove that a particular pseudocycle GW-invariant is equal to the corresponding modified polyfold GW-invariant.
All that remains is to unify these differing polyfold GW-invariants and show that they are all equal to the polyfold GW-invariant introduced at the beginning.

To begin, let us recall the definitions of the various GW-polyfolds we have encountered thus far.

\begin{enumerate}[label = (\Roman*)]
	\item \label{GW_polyfold_I}
	      The GW-polyfold of stable curves which map to the symplectic manifold $M$, and which are of class $A$, genus $0$, and with $k$ marked points:
	      \begin{gather*}
		      \cZ_{A,0,k} =
		      \left\{
		      \begin{array}{c}
			      u: S^2 \to M \\
			      \{z_1,\ldots,z_k\} \in S^2
		      \end{array}
		      \middle|
		      \begin{array}{c}
			      u \in W^{2,3}(S^2,M), \\ u_* [S^2] = A,\\
			      z_i \neq z_j \text{ if } i \neq j
		      \end{array}
		      \right\}
		      \biggm/
		      \begin{array}{c}
			      \phi \in \text{PSL}(2,\C) \\
			      \phi(z_i) = z'_i
		      \end{array}
		      \\
		      \bigsqcup \ \{\text{stable nodal curves}\}.
	      \end{gather*}
	      The perturbed GW-moduli space is given by
	      \[
		      \CM_{A,0,k}(p) = (\delbarj +p)^{-1}(0) \subset \cZ_{A,0,k}.
	      \]
	      The GW-invariants \eqref{eq:full-polyfold-gw-invariant} are defined as either
	      \begin{align*}
		       & \text{polyfold-}\GW_{A,0,k} ([\cX_1],\ldots,[\cX_k]) \\
		       & \qquad =
		      ev({\CM_{A,0,k}(p)}) \cdot \left(\cX_1 \times\cdots\times \cX_k \right).
	      \end{align*}
	      or with contributions from homology classes in $H_*(\dmspace_{0,k};\Q)$ as
	      \begin{align*}
		       & \text{polyfold-}\GW_{A,0,k} ([\cX_1],\ldots,[\cX_k];[\cB])                                               \\
		       & \qquad	=	(ev\times \pi )({\CM_{A,0,k}(p)}) \cdot \left(\cX_1 \times\cdots\times \cX_k \times \cB\right).
	      \end{align*}

	\item \label{GW_polyfold_II}
	      The GW-polyfold of graphs with three fixed marked points, i.e., the GW-polyfold of stable curves which map to the symplectic manifold $\widetilde{M} = S^2 \times M$, and which are of class $\widetilde{A} = [S^2 \times \text{pt}] + \iota_* A]$, genus $0$, and with $k$ marked points, where in addition we fix the first three marked points at $0,1,\infty$ via \eqref{eq:GW-polyfold_with_fixed_marked_points}:
	      \begin{gather*}
		      \cZ_{\widetilde{A},0,k}^{(0,1,\infty)} =
		      \left\{
		      \begin{array}{c}
			      \tilde{u}: S^2 \to \widetilde{M} \\
			      \{z_1,\ldots,z_k\} \in S^2
		      \end{array}
		      \middle|
		      \begin{array}{c}
			      \tilde{u} \in W^{2,3}(S^2,\widetilde{M}), \\
			      \tilde{u}_* [S^2] = \widetilde{A},        \\
			      z_i \neq z_j \text{ if } i \neq j,        \\
			      z_1 = 0,\ z_2 = 1,\ z_3 = \infty
		      \end{array}
		      \right\}
		      \biggm/
		      \begin{array}{c}
			      \phi \in \text{PSL}(2,\C) \\
			      \phi(z_i) = z'_i
		      \end{array}
		      \\
		      \bigsqcup \ \{\text{stable nodal curves}\}.
	      \end{gather*}
	      The perturbed GW-moduli space is given by
	      \[
		      \CM_{\widetilde{A},0,k}^{(0,1,\infty)}(p) = \big(\delbar_{\widetilde{J}}|_{\cZ_{\widetilde{A},0,k}^{(0,1,\infty)}} + p \big)^{-1} (0) \subset \cZ_{\widetilde{A},0,k}^{(0,1,\infty)}
	      \]
	      for an almost complex structure of the form $\widetilde{J} = i \oplus J_z$. The GW-invariants are defined by \eqref{eq:polyfold-graph-gw-invariants} as
	      \begin{align*}
		       & \text{polyfold-graph-}\GW_{A,0,k} ([\cX_1],\ldots,[\cX_k])                                                              \\
		       & \qquad	=	ev\left({\CM_{\widetilde{A},0,k}^{(0,1,\infty)}(p)}\right) \cdot \left(\cX_1 \times\cdots\times \cX_k \right).
	      \end{align*}
	      or by \eqref{eq:full-polyfold-graph-gw-invariants} with contributions from homology classes in $H_*(\dmspace_{0,k};\Q)$, as
	      \begin{align*}
		       & \text{polyfold-graph-}\GW_{A,0,k} ([\cX_1],\ldots,[\cX_k]; [\cB])                                                                              \\
		       & \qquad	=	(ev\times \pi)\left({\CM_{\widetilde{A},0,k}^{(0,1,\infty)}(p)}\right) \cdot \left(\cX_1 \times\cdots\times \cX_k \times \cB \right).
	      \end{align*}

	\item \label{GW_polyfold_III}
	      The GW-polyfold of graphs with three fixed marked points and which is restricted to a fixed tree, i.e., the GW-polyfold of stable curves which map to the symplectic manifold $\widetilde{M}$, and which are of class $\widetilde{A}$, genus $0$, and with $k$ marked points, where in addition we fix the first three marked points at $0,1,\infty$ via \eqref{eq:GW-polyfold_with_fixed_marked_points}, and where we restrict to a fixed tree $T$ via \eqref{eq:GW-polyfold_restricted_to_tree}:
	      \begin{gather*}
		      \cZ_{\widetilde{A},0,k}^{(0,1,\infty)} \big|_T=
		      \left\{
		      \begin{array}{c}
			      \tilde{u}: S^2 \to \widetilde{M} \\
			      \{z_1,\ldots,z_k\} \in S^2
		      \end{array}
		      \hspace{-2pt}
		      \middle|
		      \hspace{-2pt}
		      \begin{array}{c}
			      \tilde{u} \in W^{2,3}(S^2,\widetilde{M}), \\
			      \tilde{u}_* [S^2] = \widetilde{A},        \\
			      z_i \neq z_j \text{ if } i \neq j,        \\
			      z_1 = 0,\ z_2 = 1,\ z_3 = \infty,         \\
			      \pi([\tilde{u},z_1,\ldots,z_k]) \in \CM_{0,T}
		      \end{array}
		      \right\}
		      \hspace{-1pt}
		      \biggm/
		      \hspace{-4pt}
		      \begin{array}{c}
			      \phi \in \text{PSL}(2,\C) \\
			      \phi(z_i) = z'_i
		      \end{array}
		      \\
		      \bigsqcup \ \{\text{stable nodal curves}\}.
	      \end{gather*}
	      It is a subpolyfold of $\cZ_{\widetilde{A},0,k}^{(0,1,\infty)}$. The perturbed GW-moduli space is given by
	      \[
		      \CM_{\widetilde{A},0,k}^{(0,1,\infty)}\big|_T(p) = \bigm(\delbar_{\widetilde{J}}|_{\cZ_{\widetilde{A},0,k}^{(0,1,\infty)}\big|_T} + p \bigm)^{-1} (0) \subset \cZ_{\widetilde{A},0,k}^{(0,1,\infty)}\big|_T.
	      \]
	      By considering different trees and with associated GW-polyfolds of this type, the GW-invariants are defined by \eqref{eq:polyfold-graph-tree-gw-invariants} as
	      \begin{align*}
		       & \text{polyfold-fixed-tree-}\GW_{A,0,k} ([\cX_1],\ldots,[\cX_k];[\CM_{0,T}])                                                      \\
		       & \qquad	:=	ev\left({\CM_{\widetilde{A},0,k}^{(0,1,\infty)}\big|_T (p)}\right) \cdot \left(\cX_1 \times\cdots\times \cX_k \right).
	      \end{align*}
\end{enumerate}

We prove that all of these polyfold GW-invariants are equal to the original polyfold GW-invariant \eqref{eq:full-polyfold-gw-invariant}, here corresponding to \ref{GW_polyfold_I}.
In \S\ref{subsec:naturality_of_polyfold_GW-invariants} we prove that the above invariants \ref{GW_polyfold_I} are equal to the invariants \ref{GW_polyfold_II}, and in \S\ref{subsec:submanifold_GK-space_versus_restriction_to_fixed_tree} we prove that invariants \ref{GW_polyfold_II} are equal to the invariants \ref{GW_polyfold_III}.

\subsection{Naturality of polyfold invariants for the Gromov--Witten polyfold of graphs}
\label{subsec:naturality_of_polyfold_GW-invariants}

We prove that the polyfold GW-invariants \ref{GW_polyfold_I} are equal to the polyfold graph GW-invariants \ref{GW_polyfold_II}.
We have defined two distinct polyfolds which model the GW-moduli space of stable curves in slightly different ways.
To prove they are equal, we need to invoke the abstract methods of \cite{schmaltz2019naturality} concerned with showing that polyfold invariants are	\emph{natural}, and do not depend on how a moduli space is modeled.

\subsubsection{A natural inclusion map between polyfolds}
First, we need to define a natural inclusion map
\[
	I : \cZ_{A,0,k} \hookrightarrow \cZ_{\widetilde{A},0,k}^{(0,1,\infty)}
\]
which extends the map $\Psi: \cM^*_{A,0,k}(J) \to \cM_{\widetilde{A},0,k}^{*(0,1,\infty)} (i\oplus J)$ considered in \S\ref{subsec:domain-dependent-perturbation}.
We begin by giving alternative descriptions of the GW-polyfolds from that of \cite{HWZGW} to place a greater emphasis on the tree structure of a nodal curve.
To do this, we decompose the GW-polyfold $\cZ_{A,0,k}$ into a disjoint union over GW-stable labeled trees,
\[
	\cZ_{A,0,k} = \bigsqcup_T \cZ_T.
\]
Consider a tree $T$ whose vertices are labeled by a decomposition of the homology class $A$ and by the number of marked points $k$, i.e.,
\[
	A = \sum_{\alpha\in T} A_\alpha, \qquad k = \sum_{\alpha\in T} k_\alpha
\]
where $\langle [\ww], A_\alpha \rangle \geq 0$ and $k_\alpha \geq 0$ for all $\alpha \in T$.
We say such a labeled tree is \textbf{GW-stable} if it satisfies
\[
	k_\alpha + \# \{\beta \mid \alpha E \beta\} \geq 3 \quad \text{or} \quad \langle [\ww], A_\alpha \rangle \neq 0 \quad \text{for all}\quad \alpha \in T.
\]
Given such a tree we now consider $\cZ_T$, the set stable curves and which is restricted to a fixed tree $T$ \emph{without degeneration}, defined as follows:
\begin{gather*}
	\cZ_T :=
	\left\{
	\begin{array}{c}
		u_\alpha: S^2 \to M                 \\
		\{z_1,\ldots,z_{k_\alpha}\} \in S^2 \\
		\{z_\beta\}_{\alpha E \beta} \in S^2
	\end{array}
	\middle|
	\begin{array}{c}
		u_\alpha \in W^{2,3}(S^2,M),             \\ (u_\alpha)_* [S^2] = A_\alpha,\\
		u_\alpha (z_\beta) = u_\beta (z_\alpha), \\
		\text{marked, nodal pts distinct}
	\end{array}
	\right\}_{\alpha \in T}
	\biggm/
	G_T,
\end{gather*}
where $G_T$ consists of tuples $(f, \{\phi_\alpha\}_{\alpha \in T})$ where $f:T\to T$ is a tree isomorphism and $\phi_\alpha :S^2 \to S^2$ are M\"obius transformations which both preserve the marked points (with orderings) and nodal points.
Notice that the above definition of $\cZ_T$ does not contain further degeneration of a stable curve into nodal stable curves which are not modeled on the tree $T$.

We may similarly decompose the GW-polyfold of graphs with three fixed marked points $\cZ_{\widetilde{A},0,k}^{(0,1,\infty)}$.
To begin, observe that any decomposition of $\widetilde{A} = [S^2 \times \text{pt}] + \iota_* A$ must be of the form
\[
	\widetilde{A} = [S^2\times \text{pt}] + \sum_{\alpha} \iota_* A_\alpha, \quad \text{where} \quad \langle [\ww], A_\alpha \rangle \geq 0,
\]
thus we consider a GW-stable labeled tree $T$ with a distinguished vertex $\Delta \in T$ labeled by $[S^2 \times \text{pt}]$ and $k_\Delta$.
Next, we distinguish by $\alpha_0,\alpha_1,\alpha_\infty \in T$ the vertices which contain the first three marked points, which we denote as $z_0,z_1,z_\infty$. (These vertices need not be distinct.)
Recall that we fix these marked points by means of the map $\Theta$; thus we must have
\begin{equation}\label{eq:constraint}
	p_1 (\tilde{u}_{\alpha_0} (z_0)) = 0, \quad p_1(\tilde{u}_{\alpha_1}(z_1)) = 1, \quad  p_1(\tilde{u}_{\alpha_\infty}(z_\infty)) = \infty.
\end{equation}
With these considerations in mind, we may write $\cZ_{\widetilde{A},0,k}^{(0,1,\infty)} = \bigsqcup_T \cZ^{(0,1,\infty)}_T$ where
\begin{align*}
	\cZ^{(0,1,\infty)}_T :=
	 & \left\{
	\begin{array}{c}
		\tilde{u}_\Delta: S^2 \to \widetilde{M} \\
		\{z_1,\ldots,z_{k_\Delta}\} \in S^2     \\
		\{z_\beta\}_{\Delta E \beta} \in S^2
	\end{array}
	\middle|
	\begin{array}{c}
		\tilde{u}_\alpha \in W^{2,3}(S^2,\widetilde{M}),         \\
		(\tilde{u}_\Delta)_*[S^2] = [S^2\times \text{pt}],       \\
		\tilde{u}_\Delta (z_\beta) = \tilde{u}_\beta (z_\Delta), \\
		\text{marked, nodal pts distinct},                       \\
		\eqref{eq:constraint} \text{ if } \Delta = \alpha_0,\ \alpha_1, \text{ or } \alpha_\infty
	\end{array}
	\right\}
	\bigsqcup
	\\
	 & \left\{
	\begin{array}{c}
		\tilde{u}_\alpha: S^2 \to \widetilde{M} \\
		\{z_1,\ldots,z_{k_\alpha}\} \in S^2     \\
		\{z_\beta\}_{\alpha E \beta} \in S^2
	\end{array}
	\middle|
	\begin{array}{c}
		\tilde{u}_\alpha \in W^{2,3}(S^2,\widetilde{M}),         \\
		(\tilde{u}_\alpha)_*[S^2] = \iota_* A_\alpha,            \\
		\tilde{u}_\alpha (z_\beta) = \tilde{u}_\beta (z_\alpha), \\
		\text{marked, nodal pts distinct},                       \\
		\eqref{eq:constraint} \text{ if } \alpha = \alpha_0,\ \alpha_1, \text{ or } \alpha_\infty
	\end{array}
	\right\}_{\alpha \in T\setminus \Delta}
	\hspace{-20pt}
	\biggm/
	G_T,
\end{align*}
where $G_T$ consists of tuples $(f, \{\phi_\alpha\}_{\alpha \in T})$ where $f:T\to T$ is a tree isomorphism (which necessarily fixes the vertex $\Delta$) and $\phi_\alpha :S^2 \to S^2$ are M\"obius transformations which preserve the marked points (with orderings) and nodal points.

Having presented these decompositions, we may now turn to the problem of how to define the desired natural inclusion map.
We assert that any such inclusion map should restrict to our above decompositions, i.e., we should have a well-defined restriction $I |_{\cZ_T} : \cZ_T \hookrightarrow \cZ_T^{(0,1,\infty)}$.
To define such a restriction, we must answer the following questions:
\begin{itemize}
	\item Given a stable curve in $\cZ_T$, what node should map to the node distinguished by the vertex $\Delta \in T$ in $\cZ_T^{(0,1,\infty)}$?
	\item Given a collection of maps $u_\alpha: S^2 \to M$, $\alpha \in T$ with nodal constraints, how should one associate a collection of maps $\tilde{u}_{\beta} : S^2 \to \widetilde{M}$, $\beta \in T$ with nodal constraints which moreover fixes the first three marked points?
\end{itemize}

\subsubsection{Aside on basic graph theory}\label{subsubsec:graph-theory}
Let $T$ be a tree. Given any two vertices $\alpha, \beta \in T$ there exists a unique path between these vertices, which we denote by $[\alpha,\beta]$. The following observation is immediate.
\begin{proposition}{\cite{MSbook}*{Exer.~D.2.4}}
	\label{prop:basic_graph_theory}
	Given three (not necessarily distinct) vertices on a tree, $\alpha, \beta, \gamma \in T$, the intersection of connecting paths $[\alpha,\beta] \cap [\beta,\gamma] \cap [\alpha,\gamma]$ consists of a single vertex, which we denote by $\Delta\in T$.
\end{proposition}

\begin{figure}[ht]
	\centering
	\includegraphics{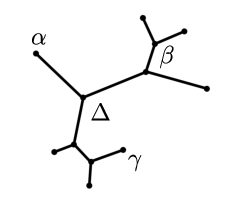}
	\caption{The distinguished vertex $\Delta$}
\end{figure}

We also introduce the following definition in order to decompose a tree. Given a tree $T$, fix a vertex $\delta \in T$. For every adjacent vertex $\alpha \in T$ which satisfies $\delta E \alpha$, we associate the subtree $T_\alpha$ defined by the union of all paths which start at $\delta$ and which pass through $\alpha$, i.e.,
\[
	T_\alpha := \bigcup \{[\delta, \gamma] \mid \gamma \in T, \ \alpha \in [\delta, \gamma] \}.
\]
Therefore having fixed a vertex $\delta$, we may then decompose the tree $T$ as follows:
\[
	T = \bigcup_{\{\alpha \mid \delta E\alpha \}} T_\alpha.
\]
This decomposition satisfies $T_\alpha \cap T_{\alpha'} = \delta$ for all $\alpha,\alpha'\in T$ such that $\delta E\alpha,$ $\delta E \alpha',$ and $\alpha \neq \alpha'.$

\subsubsection{The inclusion map on the top stratum}

The top stratum of both GW-polyfolds corresponds to the trivial tree consisting of a single vertex, $T = \{ * \}$.
We define the inclusion map $I|_{\cZ_{*}} : \cZ_{*} \hookrightarrow \cZ_{*}^{(0,1,\infty)}$ by
\[ I([(u,z_1,\ldots,z_k)]) = [\left((\id, u\circ \phi^{-1} ), \phi(z_1),\ldots, \phi(z_k) \right)]	\]
where $\phi :S^2 \to S^2$ is the unique M\"obius transformation such that $\phi(z_1) = 0,\ \phi(z_2) = 1,\ \phi(z_3) = \infty$.
We can view this as the natural extension of the diffeomorphism \eqref{eq:inclusion-map-top-stratum-gw-moduli} considered on the unperturbed GW-moduli spaces.

\subsubsection{The inclusion map on a fixed tree}

We now define $I$ for a fixed tree in the decomposition, yielding a well-defined restriction
\[
	I|_{\cZ_T}: \cZ_T \hookrightarrow \cZ_T^{(0,1,\infty)}.
\]
Denote by $\alpha_0,\alpha_1,\alpha_\infty \in T$ the vertices which contain the first three marked points.
By Prop.~\ref{prop:basic_graph_theory} may associate a vertex $\Delta\in T$ such that  $[\alpha_0,\alpha_1] \cap [\alpha_0,\alpha_\infty] \cap [\alpha_1,\alpha_\infty] = \Delta$.
For every marked point, there are two possibilities:
\begin{enumerate}
	\item the marked point lies on the component $\Delta$, or
	\item there is a unique vertex adjacent to $\Delta$ which lies on the unique path from $\Delta$ to the marked point.
\end{enumerate}
Thus, for the first three marked points we may denote by $z_{\beta_0}, z_{\beta_1}, z_{\beta_\infty} \in S^2$ for either the marked point itself, or the nodal point associated to the vertex adjacent to $\Delta$ on the path from $\Delta$ to the marked point.

We define $I$ on the distinguished vertex $\Delta$ by:
\[
	\left\{
	\begin{array}{c}
		u_\Delta: S^2 \to M                 \\
		\{z_1,\ldots,z_{k_\Delta}\} \in S^2 \\
		\{z_\beta\}_{\Delta E\beta}\in S^2
	\end{array}
	\right\}
	\mapsto
	\left\{
	\begin{array}{c}
		(\id, u_\Delta\circ \phi^{-1}) :S^2 \to S^2\times M \\
		\{\phi(z_1),\ldots,\phi(z_{k_\Delta})\} \in S^2     \\
		\{\phi(z_\beta)\}_{\Delta E\beta}\in S^2
	\end{array}
	\right\}
\]
where $\phi :S^2 \to S^2$ is again the unique M\"obius transformation such that $\{z_{\beta_0}, z_{\beta_1}, z_{\beta_\infty} \} \mapsto \{0,1,\infty\}$.
If we decompose the tree as $T = \cup_{\{ \gamma \mid \Delta E \gamma \}} T_\gamma$, then any other vertex $\alpha$ must lie on a subtree $T_\gamma$ for some vertex $\gamma$ adjacent to $\Delta$. We define $I$ on the vertex $\alpha$ by:
\[
	\left\{
	\begin{array}{c}
		u_\alpha: S^2 \to M                 \\
		\{z_1,\ldots,z_{k_\alpha}\} \in S^2 \\
		\{z_\beta\}_{\alpha E\beta}\in S^2
	\end{array}
	\right\}
	\mapsto
	\left\{
	\begin{array}{c}
		(\phi(z_\gamma), u_\alpha) :S^2 \to S^2\times M \\
		\{z_1,\ldots,z_{k_\alpha}\} \in S^2             \\
		\{z_\beta\}_{\alpha E\beta}\in S^2
	\end{array}
	\right\}.
\]
Here $z_\gamma \in S^2$ is the nodal point associated to the vertex $\Delta$ coming from the edge relation $\Delta E \gamma$, thus satisfying $u_\Delta(z_\gamma) = u_\gamma(z_\Delta)$.

Some explanation is in order. For every vertex $\alpha$, $\alpha \neq \Delta$ on the subtree $T_\gamma$, given a map $u_\alpha$ we associate a map to the fiber $\{\phi(z_\gamma)\} \times M \subset S^2 \times M$. As a consequence, the nodal constraints are satisfied. To see this, observe that given two distinct vertices $\alpha, \beta \in T$ with $\alpha E \beta$ there are two possibilities:
\begin{enumerate}
	\item $\alpha, \beta \neq \Delta$. Then $\alpha,\ \beta$ belong to the same subtree $T_\gamma$ for some $\gamma \in T$.
	\item $\alpha = \Delta$ or $\beta = \Delta$.
\end{enumerate}
From this observe the nodal constraints in the domain imply the nodal constraints in the image:
\begin{enumerate}
	\item $u_\alpha(z_\beta) = u_\beta(z_\alpha)$ implies $(\phi(z_\gamma), u_\alpha(z_\beta)) = (\phi(z_\gamma), u_\beta(z_\alpha))$,
	\item $u_\Delta(z_\alpha) = u_\alpha(z_\Delta)$ implies $(\id, u_\Delta \circ \phi^{-1}) (\phi(z_\alpha)) = (\phi(z_\alpha), u_\alpha (z_\Delta))$.
\end{enumerate}
Finally note that by the definition of $\phi$, the subtrees which contain the first three marked points, $T_{\beta_0}, T_{\beta_1}, T_{\beta_\infty}$, all map to the fibers $\{0\} \times M, \{1\} \times M, \{\infty\} \times M$. Hence observe that the constraint \eqref{eq:constraint} is satisfied.

We have only described the inclusion map on the underlying sets of the GW-polyfolds. We assert that this map is in fact $\ssc$-smooth; this can be explicitly checked according to the usual description of a local uniformizer for the GW-polyfolds as in \cite{HWZGW}.

\begin{figure}[ht]
	\centering
	\includegraphics[scale=0.8]{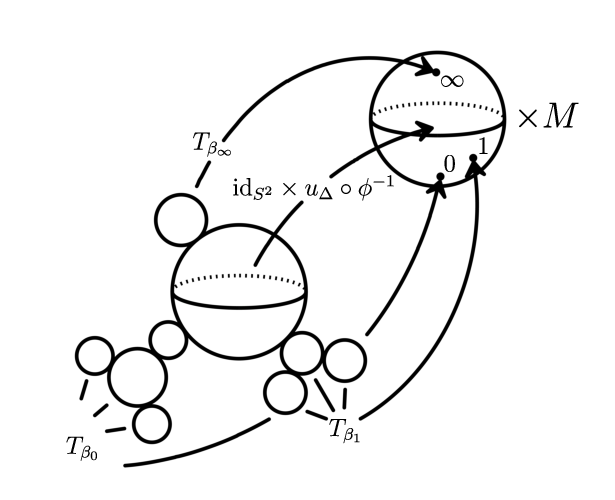}
	\caption{The inclusion map on a fixed tree}
\end{figure}

\subsubsection{Commutative diagrams of maps between polyfolds and naturality of polyfold invariants}

Having defined the above inclusion map, we may consider the following commutative diagrams.
\[\begin{tikzcd}
		\cW_{A,0,k} \arrow[r, hook] \arrow[d, "\delbarj \quad"'] & \cW_{\widetilde{A},0,k}^{(0,1,\infty)} \arrow[d, "\quad \delbar_{i\oplus J}"] &  \\
		\cZ_{A,0,k} \arrow[r, hook, "I"] \arrow[u, bend left] & \cZ_{\widetilde{A},0,k}^{(0,1,\infty)}
		\arrow[u, bend right] &
	\end{tikzcd}\]
in addition to a commutative diagram:
\[\begin{tikzcd}
		&  & M^k \times \dmlog_{0,k} \\
		\cZ_{A,0,k} \arrow[r, hook, "I"'] \arrow[rru, "ev\times\pi"] & \cZ_{\widetilde{A},0,k}^{(0,1,\infty)}
		\arrow[ru, "ev\times\pi"'] &
	\end{tikzcd}
\]
One should immediately note that $\cW_{A,0,k}$ is not the pullback bundle of $\cW_{\widetilde{A},0,k}^{(0,1,\infty)}$; it does not contain any anti-linear sections in the directions of the component $S^2$ of $S^2\times M$.

The critical hypothesis for asserting the naturality of polyfold invariants in \cite{schmaltz2019naturality} is the existence of a ``intermediary subbundle'' \cite{schmaltz2019naturality}*{Def.~3.15}.

\begin{proposition}
	The strong polyfold subbundle $\cW_{TM} \subset \cW_{\widetilde{A},0,k}$ is an intermediary subbundle.
\end{proposition}
\begin{proof}
	Given a map $\tilde{u}:S^2 \to S^2\times M$ we have a splitting of complex vector bundles $\tilde{u}^*{T(S^2\times M)} = \tilde{u}^*TS^2 \oplus \tilde{u}^* TM$.
	Likewise, given a stable curve $[\tilde{u}_\alpha]\in \cZ_{\widetilde{A},0,k}$ we have a splitting of the strong polyfold bundle $\cW_{\widetilde{A},0,k} = \cW_{TS^2}\oplus \cW_{TM}$. We may decompose the fiber over the top stratum of non-noded curves as:
	\[
		(\cW_{\widetilde{A},0,k})_{[\tilde{u}]} = W^{2,2}(S^2, \Lambda^{0,1}\otimes_i \tilde{u}^* TS^2) \oplus W^{2,2}(S^2, \Lambda^{0,1}\otimes_J \tilde{u}^* TM) / \text{PSL}(2;\C).
	\]
	The linearization of $\delbar_{i\oplus J}$ is surjective onto the first factor.
	We can therefore find a collection of vectors within $\cW_{TM}$ which span the cokernel of the linearization of $\delbar_{i\oplus J}$.
	For details on this, see furthermore \cite{MSbook}*{pg.~173}.
\end{proof}

We may therefore apply \cite{schmaltz2019naturality}*{Thm.~1.3}, and obtain the following theorem implying the polyfold GW-invariants \ref{GW_polyfold_I} and \ref{GW_polyfold_II} are equal.

\begin{theorem}
	\label{thm:GW_I_equals_GW_II}
	The polyfold GW-invariants associated to the GW-polyfolds $\cZ_{A,0,k}$ and $\cZ_{\widetilde{A},0,k}^{(0,1,\infty)}$ are equal with or without contributions from homology classes in the GK-space, i.e.,
	\[
		\text{polyfold-}\GW_{A,0,k}([\cX_1],\ldots,[\cX_k]) = \text{polyfold-graph-}\GW_{A,0,k}([\cX_1],\ldots,[\cX_k]),
	\]
	or given $[\cB] \in H_*(\dmspace_{0,k};\Q)$,
	\[
		\text{polyfold-}\GW_{A,0,k}([\cX_1],\ldots,[\cX_k],[\cB]) = \text{polyfold-graph-}\GW_{A,0,k}([\cX_1],\ldots,[\cX_k],[\cB]).
	\]
\end{theorem}

\subsection{Intersection number with a submanifold in the Grothendieck--Knudsen space versus restriction to a fixed tree type}
\label{subsec:submanifold_GK-space_versus_restriction_to_fixed_tree}

We prove that the polyfold graph GW-invariants \ref{GW_polyfold_II} with contributions from a homology class in the GK-space are equal to the polyfold graph GW-invariants \ref{GW_polyfold_III} restricted to a corresponding fixed tree type.

Let $T$ be a stable $k$-labeled tree. To this tree, we associate the submanifold $\CM_{0,T} \subset \dmspace_{0,k}$, in addition to the GW-polyfold of graphs restricted to the tree, $\cZ_{\widetilde{A},0,k}^{(0,1,\infty)} \big|_T \subset \cZ_{\widetilde{A},0,k}^{(0,1,\infty)}$.

We have a natural inclusion map between GW-polyfolds:
\[
	\iota: 	\cZ_{\widetilde{A},0,k}^{(0,1,\infty)} \big|_T \hookrightarrow \cZ_{\widetilde{A},0,k}^{(0,1,\infty)}.
\]
Via \cite{schmaltz2019naturality}*{Thm.~1.7} we may pullback an abstract perturbation and thus obtain a well-defined restriction between perturbed GW-moduli spaces:
\[
	\iota: \CM_{\widetilde{A},0,k}^{(0,1,\infty)}\big|_T (\iota^* p) \hookrightarrow \CM_{\widetilde{A},0,k}^{(0,1,\infty)}(p).
\]
Moreover, $\iota$ is a homeomorphism onto $\pi^{-1}(\CM_{0,T}) \subset \CM_{\widetilde{A},0,k}^{(0,1,\infty)}(p)$.
Given a collection of submanifolds $\cX_i \subset M$ we may perturb them so that they are transverse to both of these GW-moduli spaces, i.e,
\[ev\left( \CM_{\widetilde{A},0,k}^{(0,1,\infty)}\big|_T (\iota^* p) \right) \pitchfork (\cX_1\times \cdots \times \cX_k),
	\qquad
	ev\left( \CM_{\widetilde{A},0,k}^{(0,1,\infty)}(p) \right) \pitchfork (\cX_1\times \cdots \times \cX_k). \]
We may observe that the map $\iota$ then gives a bijection between the following points of intersection:
\[
	\iota:ev^{-1}(\cX_1\times \cdots \times \cX_k) \xrightarrow{\sim} (ev\times \pi)^{-1}(\cX_1\times \cdots \times \cX_k\times \CM_{0,T}).
\]
If we assume that $\CM_{\widetilde{A},0,k}^{(0,1,\infty)}\big|_T (\iota^* p)$ and $\cX_1\times \cdots \times \cX_k$ are of complementary dimension in $M^k$ then due to transversality the above map is a bijection between finite sets.

Here we must make note of an important subtlety. Due to the fact that $\pi$ is not a submersion, we cannot achieve transversality of the map $ev\times \pi : \CM_{\widetilde{A},0,k}^{(0,1,\infty)}(p) \to M^k\times \dmspace_{0,k}$ with the product submanifold $\cX_1\times \cdots \times \cX_k \times \CM_{0,T} \subset M^k\times \dmspace_{0,k}$ without additional perturbation of the submanifold $\CM_{0,T} \subset \dmspace_{0,k}$. However if we perturb $\CM_{0,T}$ we lose the identification of the above sets of intersection.

However, even though this intersection is not transverse, it consists of a finite set of points which counted appropriately should still give the correct intersection count. It is simple enough to make this intuition rigorous; use the equivalent definition of a polyfold GW-invariant as a branched integral; comparing the branched integrals of the Poincar\'e duals of the submanifolds at the finite intersection points via the bijection $\iota$ we see that:
\begin{align*}
	 & \text{polyfold-fixed-tree-}\GW_{A,0,k}([\cX_1],\ldots,[\cX_k];[\CM_{0,T}])                                                                  \\
	 & \qquad = 	ev\left({\CM_{\widetilde{A},0,k}^{(0,1,\infty)}\big|_T (\iota^* p)}\right) \cdot \left(\cX_1 \times\cdots\times \cX_k \right)     \\
	 & \qquad = \int_{\CM_{\widetilde{A},0,k}^{(0,1,\infty)}\big|_T (\iota^* p)} ev^* (\PD [\cX_1]\wedge \cdots \wedge \PD [\cX_k])                \\
	 & \qquad = \int_{{\CM_{\widetilde{A},0,k}^{(0,1,\infty)}(p)}} ev^* (\PD [\cX_1]\wedge \cdots \wedge \PD [\cX_k]) \wedge \pi^* \PD [\CM_{0,T}] \\
	 & \qquad = \text{polyfold-graph-}\GW_{A,0,k}([\cX_1],\ldots,[\cX_k];[\CM_{0,T}]).
\end{align*}
We have therefore proved the following.

\begin{theorem}
	The polyfold GW-invariants associated to the GW-polyfolds $\cZ_{\widetilde{A},0,k}^{(0,1,\infty)}$ and $\cZ_{\widetilde{A},0,k}^{(0,1,\infty)} \big|_T$ are equal, i.e.,
	\begin{align*}
		 & \text{polyfold-graph-}\GW_{A,0,k} ([\cX_1],\ldots,[\cX_k];[\CM_{0,T}])              \\
		 & \qquad= \text{polyfold-fixed-tree-}\GW_{A,0,k}([\cX_1],\ldots,[\cX_k];[\CM_{0,T}]).
	\end{align*}
\end{theorem}

\section{The mixed Gromov--Witten invariants of Ruan--Tian}
\label{sec:the_mixed_GW-invariants_of_RT}

As a final application, we show that the mixed GW-invariants introduced by Ruan--Tian in \cite{rt1995quatumcohomology} are also equal to the polyfold GW-invariants.

As usual, we suppose that $(M,\ww)$ is a semipositive symplectic manifold of dimension $2n$.
Fix a set of pairwise distinct marked points $\bm{x} := \{x_1,\ldots, x_k\} \in S^2$ for $k \geq 3$.
From the start, we will use domain dependent perturbations of the almost complex structure, as we wish to define invariants for multiple covers as described in Problem~\ref{prob:problem1}.
We can describe a GW-moduli space with fixed marked points directly, but it will prove helpful to instead fix the first $k$ marked points by taking a preimage.

Consider the map
\[
	\pi : \CM_{\widetilde{A},0,k+l} (\widetilde{J}) \to \dmspace_{0,k},
\]
which on the top stratum of non-noded stable $\widetilde{J}$-holomorphic curves forgets the map $\tilde{u}$ and the last $l$ marked points, i.e.,
\[
	[(\tilde{u}, y_1, \ldots, y_k, z_1, \ldots, z_l)] \mapsto [(y_1, \ldots, y_k)].
\]
The collection of fixed marked points $\bm{x}$ determines a point in the top stratum of the GK-space $[(x_1,\ldots,x_k)] \in \dmspace_{0,k}$ which we denote by $[\bm{x}]$.
We define the \textbf{GW-moduli space with $k$ fixed marked points and $l$ free marked points} by:
\[
	\CM_{\widetilde{A},0,k+l} (\bm{x}; \widetilde{J}) := \pi^{-1}([\bm{x}]) \subset \CM_{\widetilde{A},0,k+l} (\widetilde{J}).
\]
Via \cite{MSbook}*{Thm.~6.7.11}, there exist regular almost complex structure on $\widetilde{M}$ of the form $\widetilde{J}=i\oplus J_z$ such that the evaluation map
\[
	ev: \cM_{\widetilde{A},0,k+l} (\bm{x}; \widetilde{J}) \to M^k \times M^l
\]
is a pseudocycle of dimension $2n + 2c_1(A) + 2l - 2k$. And as expected, different choices of regular almost complex structure and of fixed marked points yield cobordant pseudocycles.

The \textbf{pseudocycle mixed GW-invariant} is then the homomorphism
\[
	\Phi_{A,0,k,l} : H_*^{\otimes (k+l)} (M;\Q) \to \Q
\]
defined by taking the intersection number of the above pseudocycle with a basis of representing submanifolds of $M$:
\begin{align*}
	 & \Phi_{A,0,k,l} ([A_1],\ldots,[A_k],[B_1],\ldots,[B_l])                                                                                              \\
	 & \qquad = ev\left(\cM_{\widetilde{A},0,k,l} (\bm{x}; \widetilde{J})\right) \cdot ( A_1 \times \cdots \times A_k \times B_1 \times \cdots \times B_l)
\end{align*}

We now describe the polyfold analog to this setup.
Consider the map
\[
	\pi : \cZ_{\widetilde{A},0,k+l} \to \dmlog_{0,k}.
\]
Although as we have already remarked this map is not a $\ssc$-submersion, it is transverse to the point $[\bm{x}]$. Therefore by \cite{filippenko2018constrained}*{Thm.~1.5} the preimage
\[
	\cZ_{\widetilde{A},0,k + l}(\bm{x}) := \pi^{-1} ([\bm{x}])
\]
has a polyfold structure. We call this polyfold the \textbf{mixed GW-polyfold}.
Applying the abstract comparison theorem~\ref{thm:comparison-thm} to the perturbed GW-moduli space associated to this polyfold we see that
\begin{align}
	 & \underbrace{ev(\cM_{\widetilde{A},0,k}(\bm{x};\widetilde{J}))}_{\text{pseudocycle}} \cdot (A_1 \times \cdots \times A_k \times B_1\times \cdots \times B_l) \label{eq:mixed_pseudocycle_polyfold} \\
	 & \qquad = \quad
	\underset{\substack{  \text{compact weighted}                                                                                                                                                        \\ \text{branched orbifold} }}{ev \big(\underbrace{{\CM_{\widetilde{A},0,k+l} (\bm{x}; p)}}\big)} \cdot (A_1 \times \cdots \times A_k \times B_1\times \cdots \times B_l). \nonumber
\end{align}

\subsection{Comparing the polyfold Gromov--Witten invariants and the polyfold mixed Gromov--Witten invariants}

Consider a map which forgets the last $l$ marked points of the GK-space, i.e.,
\[
	ft : \dmspace_{0,k+l} \to \dmspace_{0,k}.
\]
Observe that $\bm{x}$ determines a point $[(x_1,\ldots,x_k)] \in \dmspace_{0,k}$ which we denote by $[\bm{x}]$. The above map is a submersion, and hence the preimage $\CM_{0,\bm{x}} := ft^{-1}([\bm{x}])$ is a submanifold of $\dmspace_{0,k+l}$.
We will show that the polyfold mixed GW-invariants are equal to the usual polyfold GW-invariants when restricted to the homology class $[\CM_{0,\bm{x}}] \in H_*(\dmspace_{0,k+l}; \Q)$.

To show this relation, we use can a similar approach as in \S\ref{subsec:naturality_of_polyfold_GW-invariants} to define a $\ssc$-smooth inclusion
\[
	\iota : \cZ_{\widetilde{A},0,k + l}(\bm{x}) \hookrightarrow \cZ_{\widetilde{A},0,k+l}^{(0,1,\infty)}.
\]
And just as in \S\ref{subsec:naturality_of_polyfold_GW-invariants} we may decompose these GW-polyfolds into a disjoint union over trees and describe the inclusion map on each component.

Consider the composition
\begin{equation}
	\label{eq:ft_circ_pi}
	ft \circ \pi : \cZ_{\widetilde{A},0,k+l}^{(0,1,\infty)} \to \dmspace_{0,k}
\end{equation}
and observe the following: the map $\iota$ gives an identification between $\cZ_{\widetilde{A},0,k+l}(\bm{x})$ and the preimage $(ft \circ \pi)^{-1} ([\bm{x}]) =  \pi^{-1}(\CM_{0,\bm{x}}) \subset \cZ_{\widetilde{A},0,k+l}^{(0,1,\infty)}$.

Unlike in \S\ref{subsec:submanifold_GK-space_versus_restriction_to_fixed_tree}, one may observe that the map~\eqref{eq:ft_circ_pi} is in fact transverse to the point $[\bm{x}] \in \dmspace_{0,k}$, hence also $\pi\left(\cZ_{\widetilde{A},0,k+l}^{(0,1,\infty)}\right) \pitchfork \CM_{0,\bm{x}} $.
Moreover, given a collection of submanifolds $A_i, B_j \subset M$, we assert that
\begin{align*}
	 & ev\left(\cZ_{\widetilde{A},0,k + l}(\bm{x}) \right)\pitchfork (A_1\times \cdots \times A_k\times B_1\times \cdots\times B_l),                                               \\
	 & \qquad (ev\times \pi)\left(\cZ_{\widetilde{A},0,k+l}^{(0,1,\infty)}\right) \pitchfork (A_1\times \cdots \times A_k\times B_1\times \cdots\times B_l \times \CM_{0,\bm{x}}).
\end{align*}
Via \cite{schmaltz2019steenrod}*{Prop.~3.10} we therefore have the freedom to choose abstract perturbations which yield perturbed GW-moduli spaces which are will also be transversal to these submanifolds.
Thus we may pullback an abstract perturbation and obtain a well-defined restriction between perturbed GW-moduli spaces:
\[
	\iota : \CM_{\widetilde{A},0,k+l}(\bm{x};\iota^*p) \hookrightarrow \CM_{\widetilde{A},0,k+l}^{(0,1,\infty)}(p).
\]
such that moreover
\begin{align*}
	 & ev\left(\CM_{\widetilde{A},0,k+l}(\bm{x};\iota^*p)\right)\pitchfork (A_1\times \cdots \times A_k\times B_1\times \cdots\times B_l),                                            \\
	 & \qquad (ev\times \pi)\left(\CM_{\widetilde{A},0,k+l}^{(0,1,\infty)}(p)\right) \pitchfork (A_1\times \cdots \times A_k\times B_1\times \cdots\times B_l \times \CM_{0,\bm{x}}).
\end{align*}
We may use a polyfold change of variables theorem \cite{schmaltz2019naturality}*{Thm.~2.47} to conclude
\begin{align}
	 & ev\left(\CM_{\widetilde{A},0,k+l}(\bm{x};\iota^*p)\right)\cdot (A_1\times \cdots \times A_k\times B_1\times \cdots\times B_l), \label{eq:mixed_versus_graph}                          \\
	 & \qquad = (ev\times \pi)\left(\CM_{\widetilde{A},0,k+l}^{(0,1,\infty)}(p)\right) \cdot (A_1\times \cdots \times A_k\times B_1\times \cdots\times B_l \times \CM_{0,\bm{x}}). \nonumber
\end{align}
Combining equations~\eqref{eq:mixed_pseudocycle_polyfold} and \eqref{eq:mixed_versus_graph} with theorem~\ref{thm:GW_I_equals_GW_II} we conclude:

\begin{theorem}
	Fixing a collection of marked points $\bm{x}$ corresponds to pairing with the homology class $[\CM_{0,\bm{x}}] \in H_*(\dmspace_{0,{k+l}};\Q)$.
	The pseudocycle mixed GW-invariants are related to the polyfold GW-invariants via the equation:
	\begin{align*}
		 & \Phi_{A,0,k,l} ([A_1],\ldots,[A_k],[B_1],\ldots,[B_l])                                           \\
		 & \qquad= \text{polyfold-}\GW_{A,0,k+l} ([A_1],\ldots,[A_k],[B_1],\ldots,[B_l], [\CM_{0,\bm{x}}]).
	\end{align*}
\end{theorem}

\section*{Acknowledgments}
The author thanks Benjamin Filippenko for discussions explaining his work on a general transverse preimage theorem for polyfold theory.
The author also thanks {\"O}zlem Y{\"o}nder for her help with the illustrations.



\end{document}